\RequirePackage{fix-cm}

\documentclass[onefignum,onetabnum]{siamart220329}

\usepackage{microtype}
\usepackage{xurl}

\usepackage{tikz}%
\usetikzlibrary{plotmarks}
\usepackage{caption}%for captions inside minipage i.e. Nodes and Element table
\usepackage{epstopdf}

\usepackage[colorinlistoftodos]{todonotes} %adds to do comments
\usepackage{color,xcolor}
\usepackage{amsmath,amssymb,bm}
\usepackage{graphicx,graphics}
\usepackage{subcaption} % modern replacement for subfigure/subfig
\definecolor{refgreen}{rgb}{0,0.5,0}
\usepackage{booktabs} % for tables

\usepackage{pgfplots}
\usepackage{grffile}
\pgfplotsset{compat=newest}
\usetikzlibrary{plotmarks}
\usetikzlibrary{arrows.meta}
\usepgfplotslibrary{patchplots}

\date{ \phantom{b} \today \phantom{e}}
\newcommand{\R}{\mathbb{R}}

\usepackage{hyperref}
\usepackage[normalem]{ulem}

\usepackage{algorithm} %For algorithms
\usepackage{algorithmic}

\usepackage[export]{adjustbox} % for trim and clip keys in \includegraphics

\usepackage{diffcoeff}
\usepackage{bm} %big font for math symbols

\definecolor{refblue}{rgb}{0,0,0.75}
\definecolor{refblueb}{rgb}{0,0,1}
\definecolor{refgreen}{rgb}{0.13,0.55,0.13}
\definecolor{refred}{rgb}{1,0,0}

\hypersetup{
	colorlinks   = true, %Colours links instead of ugly boxes
	urlcolor     = refblue, %Colour for external hyperlinks
	linkcolor    = refblueb, %Colour of internal links
	citecolor   = refgreen %Colour of citations
}

\usepackage[T1]{fontenc}
\usepackage{textcomp}
\usepackage{lmodern}
\usepackage{hyperref}

\usepackage[autostyle=true]{csquotes}

\def            \d          {\text{d}}

\renewcommand*  {\epsilon}  {\varepsilon}

\def            \i          {\text{i}}
\newcommand     {\id}       {\mathop{\text{Id}}}

\newcommand     {\inv}      {^{-1}}
\newcommand     \laplace    {\varDelta}

\newcommand     \nb         {\nabla}
\newcommand     \Om         {\varOmega}

\def            \to         {\rightarrow}

\newcommand     \vphi       {\varphi}

\newcommand\bfb{{\mathbf b}}

\newcommand\bfe{{\mathbf e}}
\newcommand\bff{{\mathbf f}}

\newcommand\bfn{{\mathbf n}}

\newcommand\bfu{{\mathbf u}}
\newcommand\bfv{{\mathbf v}}

\newcommand\bfx{{\mathbf x}}

\newcommand\bfA{{\mathbf A}}

\newcommand\bfH{{\mathbf H}}

\newcommand\bfM{{\mathbf M}}

\newcommand\andquad{\quad\hbox{ and }\quad}

\newcommand{\Ga}{\varGamma}

\newcommand{\nbg}{\nabla_{\Ga}}
\newcommand{\nbgh}{\nabla_{\Ga_h}}
\newcommand{\spn}{\textnormal{span}}

\newcommand\loc{_{\text{loc}}}

\newcommand{\hi}{\hat i}
\newcommand{\hj}{\hat j}
\newcommand{\hk}{\hat k}

\definecolor{txtgrey}{rgb}{0.55,0.55,0.55}

\usepackage{xspace}
\newcommand{\matlab}{\textsc{Matlab}\xspace}

\newcommand{\redoff}{\color{black}}

\usepackage{tikz}
\usetikzlibrary{matrix,calc}

\usepackage{booktabs}
\usepackage{siunitx}

\definecolor{matlabgreen}{rgb}{0,0.5,0}
\usepackage{listings}
\newsiamremark{remark}{Remark}

\title{$\ell$FEM: An efficient loop-free Matlab implementation \\ of isoparametric bulk and surface finite elements}

\headers{$\ell$FEM: loop-free bulk and surface FEM assembly}{B.~Kov\'acs and M.~Lantelme}

\author{Bal\'azs Kov\'acs$^*$ \and
Michael Lantelme\thanks{Institute of Mathematics, Paderborn University, Warburgerstr.~100, 33098 Paderborn, Germany} 
e-mail: \email{\{balazs.kovacs,lantelme\}@math.upb.de} }

\begin{document}

\maketitle

\begin{center}
	\emph{This paper is dedicated to the memory of Dominik Edelmann (1991--2025).}
\end{center}

\medskip

\begin{abstract}
	The $\ell$FEM \matlab package provides a simple, efficient, and flexible implementation of isoparametric finite elements in bulk domains and on surfaces. The finite element matrix assemblies are based on \matlab's paged operators and therefore completely loop-free. We give a short and conscious description of high-order isoparametric surface finite elements, which is then used to describe the assembly process and the implementation.
	We report on relevant numerical experiments (runtime comparisons, modifications for non-linear problems, etc.), and on additional functions, examples, and a testing unit which are all part of the $\ell$FEM package. 
\end{abstract}

{\small

MSC:
% as examples https://mathscinet.ams.org/mathscinet/msc/pdfs/classifications2020.pdf
35-04, % Software, source code, etc. for problems pertaining to partial differential equations
58J05, % Elliptic equations on manifolds, general theory [
35R01, % PDEs on manifolds
65-04, % Software, source code, etc. for problems pertaining to numerical analysis
%65D32, % Numerical quadrature and cubature formulas
65M60, % Finite element, Rayleigh-Ritz and Galerkin methods for initial value and initial-boundary value problems involving PDEs
65N30, % Finite element, Rayleigh-Ritz and Galerkin methods for boundary value problems involving PDEs
65Y05, % Parallel numerical computation
}

\section{Introduction}
We present the $\ell$FEM \matlab package providing \textit{simple}, \textit{efficient}, and \textit{flexible} implementation of isoparametric finite elements in smooth domains \emph{and} on smooth surfaces. $\ell$FEM provides \textit{loop-free} assembly for linear and quadratic isoparametric finite elements for surfaces of dimension 1 and 2, and for bulk domains of dimension 2 and 3.

The $\ell$FEM package is particularly advantageous for the spatial discretisation of problems which require multiple assemblies of finite element matrices, e.g., for non-autonomous or quasi-linear evolution equations, (evolving) surface partial differential equations and geometric surface flows \cite{DeckelnickDziukElliott_Acta_2005,DziukElliott_Acta_2013,MCF,BGN_survey,ElliottRanner2021,MCFdiff,anisotropicMCF}, coupled bulk--surface problems \cite{ElliottRanner2013,EylesKingStyles2019,StylesVanYperen2021,bulksurface_coupling}, and partial differential equations with dynamic boundary conditions  \cite{dynbc,HochbruckHippStoher,CH_fulldiscrete}, and also for adaptive finite element methods \cite{Demlow_Dziuk_2007,Bonito_Devore_Nochetto_2013,adaptive_surface_FEM}.

The largest efficiency improvements arise due to vectorization which is a common theme for fast FEM related codes for Euclidean domains in \matlab, see, e.g., the  \cite{funken2011,Rahman_Valdman2013,funken_beuter_2024,mallesham2025}.

\begin{equation*}
	\text{The $\ell$FEM package is available on the git repository: \href{http://go.upb.de/ellFEM}{\texttt{go.upb.de/ellFEM}}.}
\end{equation*}

The main development goals behind $\ell$FEM are \emph{efficiency}, \emph{simplicity}, and \emph{flexibility}. 
Inspired by the almost identical goals of Persson \& Strang for their meshing tool DistMesh \cite{DistMesh}, see \href{https://persson.berkeley.edu/distmesh/}{\texttt{https://persson.berkeley.edu/distmesh/}}.

\textit{Efficiency.} 
The package uses the classical reference element technique to assemble finite element matrices, but it is completely loop-free.
%but unlike traditional assembly processes, does not rely on loops over mesh elements, and/or loops over nodes of a quadrature rule.
%Instead we follow a vectorized approach, where the required transformations, quadratures, and other operations are performed simultaneously for all elements using large matrices and tensors. 
It is \emph{completely vectorized}, i.e., all operations, transformations, quadratures, etc.~are performed simultaneously for all elements using matrices and tensors. Due to its roots, the matrix--vector arithmetic of \matlab (``Matrix laboratory'') is extremely well-suited for this task.
This parallel, loop-free approach significantly reduces computational time and hence increases scalability. 

\textit{Simplicity.} 
As we will demonstrate later on, the codes of $\ell$FEM are close to the mathematical formulas, making the codes easy to understand and modify, while still not sacrificing efficiency for readability. 
Therefore, some parts could be further sped-up, which would lead away from mathematical formulation and lead to obscure codes.

\textit{Flexibility.} 
$\ell$FEM is intended to be useful in research prototyping, and to be accessible in education.
Later on we will demonstrate, that matrices for complicated non-linear terms can easily be assembled, see Section~\ref{sec:ex2_KLL}. 
Furthermore, for instance, the assembly functions can readily be extended to higher dimension, or to finite elements of polynomial degree $k \geq 3$. 

\textbf{Finite element libraries.} 
There are various finite element libraries providing highly optimized assembly codes. In the following, we first give an overview of existing \matlab-based FEM implementations before turning to general-purpose finite element solvers.

%[briefly describe \matlab stuff:]
%Up to our knowledge, existing FEM implementations in \matlab are either loop-free, but do not provide support for isoparametric elements, e.g., \cite{Pretorius_et_al,??}, or can handle isoparametric FEM, but use loops during the assembly, e.g., \cite{Bartels_et_al_2006,??,??}.

There exists \matlab-based FEM implementations which are not loop-free supporting both parametric in 2D and 3D \cite{Carstensen_1999} and even quadratic isoparametric \cite{Bartels_et_al_2006} elements in 2D. Later  packages introduced loop-free implementations for linear finite elements in 2D in \cite{funken2011, Cuvelier2013}, for 3D and higher in \cite{Chen_2009_ifem,Feifel_Funken_2024}, and for linear isoparametric elements in \cite{Rahman_Valdman2013}.

In this work, we adopt a framework similar to \cite{Rahman_Valdman2013}, but extend it in several important directions.
We provide a clear and general structure how to handle arbitrary order elements (illustrated for second order), and arbitrary quadrature rules. In particular we extend the analysis to general surface and bulk domains, which up to our knowledge was not yet done. Additionally, we leverage the implementation by utilising modern \matlab functions, which simplifies the translation of naive code, automatically optimizes linear algebra, and provides built-in GPU support.

Established finite element packages implementing isoparametric surface finite elements are
\href{https://www.alberta-fem.de/}{ALBERTA}, \href{https://www.dune-project.org/}{DUNE}, \href{https://www.firedrakeproject.org}{Firedrake}, \href{https://mfem.org/}{MFEM}, \href{https://ngsolve.org}{Netgen/NGSolve}, etc. These packages provide optimized assembly routines, which are written in low-level languages like C/C++. Yet using these packages to implement novel algorithms for systems coupling parabolic PDEs on evolving surfaces to geometric surface flows, e.g., \cite{MCFdiff,bulksurface_coupling,anisotropicMCF}, or numerical surgery \cite{numerical_surgery} may still prove to be a challenging task. $\ell$FEM provides a simple direct access to source codes, but retains quick execution.  

\textbf{Components.} 
The $\ell$FEM package provides:
\begin{itemize}
\item[-] numerical tests on surfaces and in bulk domains, including: 
elliptic problems, mean curvature flow, a runtime comparison to other finite element libraries, and a runtime test for the dependency on polynomial degree of the finite elements
\item[-] a testing unit, 
\item[-] a driver file/driver files to run all experiments presented in this paper,
\item[-] a basic documentation wiki.
\end{itemize}

\textbf{Outline.} 
The paper is organised as follows:
Section~\ref{sec:model_problem} introduces elliptic surface PDEs as a model problem. Using isoparametric finite elements the finite element matrix--vector formulation is derived. In Section~\ref{sec:assembly} we set up the framework to  determine mass and stiffness matrices, including quadrature, for surface PDEs in an elementwise manner. Based on the formulation  Section~\ref{sec:overview_ellFEM} presents a naive approach to implement the assembly, discuss its flaws and introduce remedies. Via vectorization based on $\matlab$-inbuilt functions we present an efficient, loop-free assembly code. Additionally we give an overview over the functionalities of the $\ell$FEM package.
Finally in Section~\ref{sec:num_ex} extensive numerical experiments show assembly runtime in general and compared to other available finite element packages. These examples provide insights into GPU computations, solving geometric PDEs, and how to deal with non-linear terms.

\section{Model problem: elliptic surface PDE}
\label{sec:model_problem}
To describe the concepts of our implementation, we use the Poisson problem on a two-dimensional closed surface as a model problem:
Let $\Ga \subset \R^3$ be a two-dimensional, closed, smooth (at least $C^2$) surface, whose principal curvatures are bounded, and let $f \in L^2(\Ga)$, $\int_{\Ga} f = 0$, be a given mean-free function. We consider the elliptic surface PDE: Find a sufficiently smooth solution $u \colon \Ga \to \R$ such that
\begin{equation}
\label{eq:model problem}
	-\laplace_\Ga u = f \qquad \textnormal{in $\Ga$} .
\end{equation}
The weak formulation of \eqref{eq:model problem} reads: Find $u \in H^1(\Ga)$, $\int_{\Ga} u = 0$, such that, for all $v \in H^1(\Ga)$, there holds
\begin{equation}
\label{eq:weak model problem}
	\int_\Ga \nbg u \cdot \nbg v = \int_\Ga f v .
\end{equation}
Well-posedness of the weak problem follows by the Lax--Milgram lemma. 
The tangential (or surface) gradient of a smooth scalar function $u\colon \Ga \to \R$ is given by, see, e.g.,~\cite[Section~2.1 and Definition~2.3]{DziukElliott_Acta_2013},
\begin{equation*}
	\nabla_{\Ga} u = \nabla \overline{u} - (\nabla \overline{u} \cdot \nu) \nu ,
\end{equation*} 
where $\nu \colon \Ga \to \R^3$ is the outer unit normal field to $\Ga$, and $\overline{u}$ is an extension of $u$ into a neighbourhood $U \subset \R^3$ around $\Ga$. The tangential gradient is independent of the extension. 
The above description uses standard notation for Sobolev spaces, surface PDEs, etc., for more details we refer to \cite{Dziuk88,DziukElliott_Acta_2013}.

\subsection{Isoparametric surface finite elements of order $k$}
The \textit{surface finite element method} of polynomial order $k$ (often abbreviated to: P$k$ surface FEM), following \cite{Dziuk88,Demlow2009}, seeks an approximate solution based on an admissible surface triangulation of polynomial order $k \geq 1$ denoted by $\Ga_h^k \neq \Ga$, see~\cite[Section~2.2]{Demlow2009}, which consist of elements being a degree-$k$ polynomial image of a reference element ($\hat{E}$), and whose nodes all lie on $\Ga$. As we will always deal with surface FEM of arbitrary polynomial degree we will often set $\Ga_h := \Ga_h^k$.

\begin{figure}[htbp]
    \centering
    \includegraphics[width=0.45\linewidth]{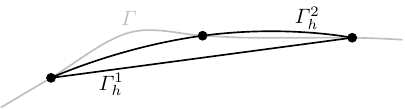}
    \hfill
    \includegraphics[width=0.45\linewidth]{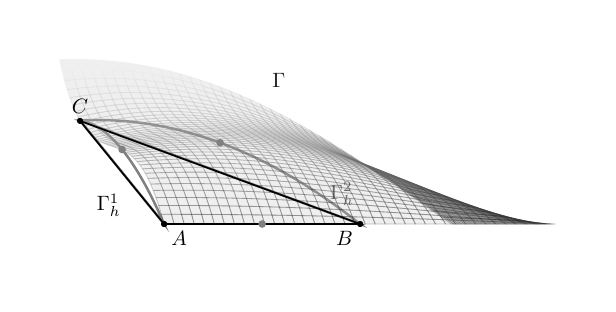}
    \caption{First and second order elements of a discrete surface $\Ga_h^k$ approximating the surface $\Ga$ in two (on the left) and three dimensions (on the right).}
    \label{fig:discrete elements}
\end{figure}

Let $S_h^k$ denote the surface finite element space on $\Ga_h$ collecting continuous functions that are elementwise of polynomial degree $k$, with Lagrangian basis
\begin{equation}
\label{eq:S_h^k basis}
	S_h^k  = \text{span} \{ \phi_1, \dotsc, \phi_{N} \} .
\end{equation}
Note that the surface FEM method is not conform, $S_h^k \nsubseteq H^1(\Ga)$.

The discrete tangential gradient is given, see, e.g.,~\cite[Section~4.4]{DziukElliott_Acta_2013}, elementwise by
\begin{equation}
\label{eq:discrete tangential gradient}
	\nabla_{\Ga_h} u_h := P_h \nabla \overline{u_h} := \nabla \overline{u_h} - (\nabla \overline{u_h} \cdot \nu_h) \nu_h  ,
\end{equation}
where $\overline{u_h}$ is an extension of $u_h \in S_h^k$ into $U$ which is constant along normals. The discrete tangential gradient is also independent of the extension. 

The degree $k$ isoparametric surface FEM discretisation then reads: Find $u_h(\cdot,t) \in S_h^k$, $\int_{\Ga_h} u_h = 0$, such that, for all $v_h \in S_h^k$,
\begin{equation}
\label{eq:discretised model problem}
	\int_{\Ga_h} \nabla_{\Ga_h} u_h \cdot \nabla_{\Ga_h} v_h = \int_{\Ga_h}  f_h v_h ,
\end{equation}
where $f_h \in S_h^k$ is a suitable approximation to $f$. 

Through the closest point projection, i.e.~via the unique solution $x^\ell \in \Ga$ to $x^\ell := x - d(x) \, \nu (x^\ell)$ for arbitrary $x \in \Ga_h \subset U$, the \emph{lift} of a function $w \colon \Ga_h \rightarrow \R$ onto $\Ga$ is given by $w^{\ell} (x^\ell) := w (x)$, while the \textit{unlift} $w^{-\ell} \colon \Ga_h \to \R$ is defined such that $(w^{-\ell})^\ell = w \colon \Ga \to \R$ holds.

High-order isoparametric error estimates for polynomial degree $k \geq 1$ were shown, see \cite[Corollary~4.2]{Demlow2009} (in slightly different form), and also \cite[Theorem~4.9]{DziukElliott_Acta_2013}: 
\begin{equation}
\label{eq:error estimates for SFEM}
	\|u - u_h^\ell\|_{L^2(\Ga)} + h \|\nbg (u - u_h^\ell)\|_{L^2(\Ga)} \leq c h^{k+1} \|f\|_{L^2(\Ga)} + c_f \, h^{k+1} ,
\end{equation}
provided that $\|f - f_h^\ell\|_{L^2(\Ga)} = c_f \, h^{k+1}$. The estimates in \cite{Demlow2009} even split the surface and finite element approximation orders.

\subsection{Matrix--vector formulation}
%The discrete problem \eqref{eq:discretised model problem} is equivalent to the linear system
%\begin{equation*}
%	\bfA \bfu = \bfb ,	
%\end{equation*}
%where $\bfu \in \R^N$ collects the nodal values of $u_h = \sum_{j=1}^{N} \bfu|_j \phi_j \in S_h^k$, and where $\bfA \in \R^{N \times N}$ is the stiffness-matrix and $\bfb \in \R^N$ is the load vector. 
%
%The mass- and stiffness-matrix, and the load vector are respectively given by
%\begin{gather*}
%	\bfM|_{ij} = \int_{\Ga_h} \phi_i \phi_j , \qquad
%	\bfA|_{ij} = \int_{\Ga_h} \nbgh \phi_i \cdot \nbgh \phi_j , 
%	\andquad 
%	\bfb|_{j} = \int_{\Ga_h} f_h \phi_j .
%\end{gather*}
%for $i,j = 1,  \dotsc,N$. Note that, since $f_h = \sum_{j=1}^{N} \bff|_j \phi_j \in S_h^k$, $\bfb = \bfM \bff$ holds.

The discrete problem is equivalent to the linear system
\begin{equation*}
	\bfA_0 \bfu = \bfb_0 ,
\end{equation*}
where $\bfu \in \R^N$ collects the nodal values of $u_h = \sum_{j=1}^{N} \bfu|_j \phi_j \in S_h^k$ with zero mean-value, and where $\bfA_0 \in \R^{N+1 \times N}$ is the stiffness-matrix and 
$\bfb_0 \in \R^{N+1}$ is the load vector, such that
\begin{equation*}
	\bfA_0 = [\bfA ; \bfe^T]
	\qquad \textnormal{and} \qquad
	\bfb_0 = [\bfb ; 0],
\end{equation*}
where $\bfe = (1,\dotsc,1)^T \in \R^N$. 
The last appended line takes care of the mean-freeness of the numerical solution.

The mass- and stiffness-matrix, and the load vector are respectively given by
\begin{gather*}
	\bfM|_{ij} = \int_{\Ga_h} \phi_i \phi_j , \qquad
	\bfA|_{ij} = \int_{\Ga_h} \nbgh \phi_i \cdot \nbgh \phi_j , \\
	\textnormal{and} \qquad
	\bfb|_{j} = \int_{\Ga_h} f_h \phi_j ,
\end{gather*}
for $i,j = 1,  \dotsc,N$. 
Note that, since $f_h = \sum_{j=1}^{N} \bff|_j \phi_j \in S_h^k$, 
$\bfb = \bfM \bff$ holds.

\section{Assembly of isoparametric surface finite element matrices}
\label{sec:assembly}
\subsection{Reference map for surface FEM}
\begin{figure}[htbp]
	\centering
	\begin{tikzpicture}
	% the 3D coordiante system
	\draw[color=gray,->] (-0.42,0) -- (4.42,0) node[pos=1]{ \large \quad $\xi_1$};	\draw[color=gray,->] (-0.125,-0.25) -- (2,4) node[pos=1]{ \large ~\quad $\xi_2$};
	\draw[color=gray,->] (0,-0.27) -- (0,4.25) node[pos=1,align=left]{\large $\xi_3$ \quad};
%	% gray curves
%	\draw[color=gray,dotted] (0,0) to[bend left] (6,1.5);
%	\draw[color=gray,dotted] (3,0) to[bend right] (9,1);
%	\draw[color=gray,dotted] (1,2) to[bend left] (7,3);
%	\draw[color=gray,dotted] (0,3) to[bend left] (6.542,5);
		% gray curves
	\draw[color=gray,dashed] (0,0) to[out=20,in=174] (7,1.5);
	\draw[color=gray,dashed] (3,0) to[out=-18,in=201] (10,1);
	\draw[color=gray,dashed] (1,2) to[out=13,in=198] (8,3);
	\draw[color=gray,dashed] (0,3) to[out=38,in=184] (7.542,5);
	% the reference triangle \hat{E}
	\draw[very thick]
	(0,0) node[align=right, above] {\large \qquad 1}
	-- (3,0) node[align=left, above] {\large 2 \qquad\,\,}
	-- (1,2) node[align=right, below] {\\ \large\ 3 }
	-- cycle;
	\draw[thick,fill=black]
	(0,0) circle (2pt) node[align=left, below] {\\ \large $(0,0,0)$}
	-- (3,0) circle (2pt) node[align=right, below] {\large $(1,0,0)$}
	-- (1,2) circle (2pt) node[align=left, above] {\large $(0,1,0)$}
	-- cycle;
	\draw
	(1.042,0.742) node{\rotatebox[origin=c]{-4}{\Large $\hat{E}$}};

	% midpoint of reference element
	\coordinate (m12) at ($(0,0)!0.5!(3,0)$);
	\coordinate (m23) at ($(3,0)!0.5!(1,2)$);
	\coordinate (m31) at ($(1,2)!0.5!(0,0)$);
	\draw[fill=black,thick] (m12) circle (2pt) node[below] {\large\ 4};
	\draw[fill=black,thick] (m23) circle (2pt) node[above] {\large\ 5};
	\draw[fill=black,thick] (m31) circle (2pt) node[above, left]  {\large\ 6};

	% the normal vector \nu_0
	\draw[very thick,->] (0,0) -- (0,3) node[pos=0.5,align=left]{\Large $\hat{\nu}$\quad\, };
	%%%
	% arrows
	\draw[->,thick] (3,1.5) .. controls +(up:0.52cm) and +(left:0.52cm) .. node[above,sloped] {$F:\hat{E} \to E$} (6.5,1.75);
	\draw[->,thick] (7.5,1.1) .. controls +(down:0.52cm) and +(right:0.52cm) .. node[below,sloped] {$F\inv:E \to \hat{E}$} (3.5,0.8);
	
	% element E
	%\draw[very thick]
	%(7,1.5) -- (10,1);
	\draw[very thick]
	(8,3) -- (7,1.5);
	\draw[thick,fill=black]
	(7,1.5) circle (2pt) node[align=left, below] {\large $x_1$ \,}
	(10,1) circle (2pt) node[align=right, below] {\, \large $x_2$}
	(8,3) circle (2pt) node[align=left, above] {\large $x_3$};
	\draw
	(8.1,2.2) node{\rotatebox[origin=c]{-4}{\Large $E$}};
	
	% midpoints of general triangle
	\coordinate (xm12) at ($(7,1.5)!0.5!(10,1) + (0.1,0.2)$);
	\coordinate (xm23) at ($(10,1)!0.5!(8,3) + (0.2,0.2)$);
	\coordinate (xm31) at ($(8,3)!0.5!(7,1.5)$);
	\draw[fill=black,thick] (xm12) circle (2pt) node[below right] {\large $x_4$};
	\draw[fill=black,thick] (xm23) circle (2pt) node[above right] {\large $x_5$};
	\draw[fill=black,thick] (xm31) circle (2pt) node[above]  {\large $x_6$};
	\draw[very thick] (10,1) .. controls ($(xm23) + (0.07,0.07)$) .. (8,3); %curved edge 235
	\draw[very thick] (10,1) .. controls ($(xm12) + (0.04,0.07)$) .. (7,1.5); %curved edge 124
	
	% arrows between midpoints
	\draw[color=gray,dashed] (m12) to[out=-30,in=220] (xm12);
	\draw[color=gray,dashed] (m23) to[out=0,in=220] (xm23);
	\draw[color=gray,dashed] (m31) to[out=0,in=180] (xm31);

	% the normal vector \nu_E
	\draw[very thick,->] (7,1.5) -- (7.542,5) node[pos=0.5,align=left]{\Large $\nu_E$\quad\, };
	\draw[thick,fill=black]
	(7.542,5) circle (2pt) node[align=right, below] {\qquad\quad\quad\, \large $x_1 + \nu_E$};
	\end{tikzpicture}
	\caption{The reference element $\hat{E}$ and an arbitrary element $E$, and the transformation maps $F$ and $F\inv$ in between.}
	\label{fig:notation}
\end{figure}
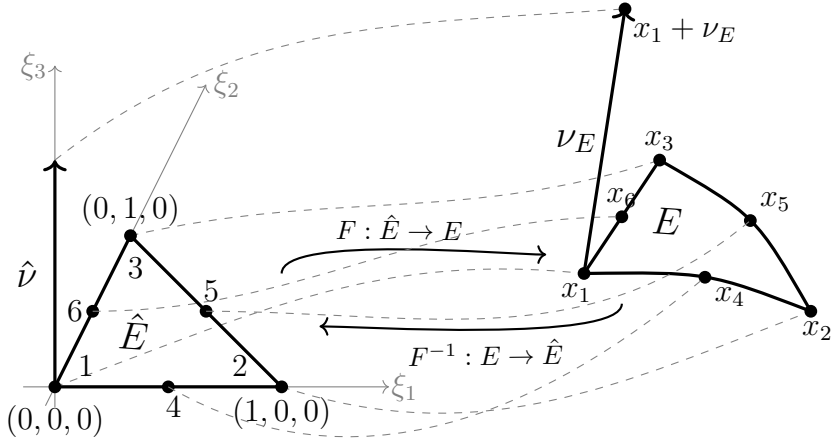

We use a reference element $\hat{E}$ with normal $\hat{\nu}$. Any element $E$ of the mesh is corresponded to the reference element by the bijective map $F_E \colon \hat E \rightarrow E$. The setting in our example is illustrated in Figure~\ref{fig:notation}. The nodes $x_j$ describe the element, we abbreviate $\bfx = (x_j)$ for the collection of nodes and employ this notational convention for further quantities. Note that for surface finite elements we have an additional degree of freedom which is used to uniquely orient the manifold. This is realised by requiring that the normals of the the reference element $\hat\nu$ and general element $\nu$ fulfil $F(\hat\nu) = x_1 + \nu_E$.

We will denote the global basis functions $\phi_j$, while on the reference element $\hat{E}$ by $\vphi_{\hj}$, where $\hj$ is the local label of the globally labelled node $j$. Recall the Lagrange basis functions satisfy $\phi_j (x_i) = \delta_{ji}$. Any element $E$, with nodes $x_j$, is also given as the image of the map 
\begin{equation}
\label{eq:F_map}
	F(x) = \sum_{\hj=1}^6 x_j \vphi_{\hj}(x) .
\end{equation}
The basis functions on the reference element, $\hat{E} = \spn\big( (1,0), (0,1), (0,0) \big) \subset \R^2$, are given by
\begin{align}
\label{eq:form_functions_s2dp2}
	\vphi_1(x,y) &= 1 - x - y, & \vphi_2(x,y) &= x, & \vphi_3(x,y) &= y, \\
	\vphi_4(x,y) &= 4x(1 - x - y), & \vphi_5(x,y) &= 4xy, & \vphi_6(x,y) &= 4y(1 - x - y).
\end{align}

\begin{remark}
	The basic concepts for the finite element assembly for bulk problems are straightforward, interpreting the domain as a flat surface with fixed normal vector $\nu = (0,\dotsc,0,1)$. 
	Therefore, we do not describe this here in detail, but refer, e.g., to \cite{Bartels_et_al_2006}.
\end{remark}

\subsection{Computing global matrices}
We use the following integral formulas to resolve the mass and stiffness matrix. For P2 elements we will require quadrature rules. We start by transforming all expression on the general element $E$ to the reference element $\hat{E}$. Any dependencies on the element are displayed as the subindex $E$.
\begin{equation}
\label{eq:basic integral formula}
	\int_{E} f(x) \d x = \int_{F_E(\hat E)} f(x) \d x = \int_{\hat E} f(F_E(\xi)) |\det DF_E| \d \xi . 
\end{equation}
Further we use the fact that the reference and global basis functions are related by $F_E$:
\begin{equation}
\label{eq:basic function formula}
	\phi_j (x) = \vphi_{\hj}({F_E}\inv(x)) : E \to \hat E \to \R .
\end{equation}

Due to the local supports of the basis function we split the mass and stiffness matrices elementwise, $j,k=1,2\dotsc,N$,
\begin{alignat}{4}
\label{eq:M matrix}
    \bfM|_{kj} &= \int_{\Ga_h} \phi_j \phi_k \d x 
    &\ = &\  \sum_E \int_E \phi_j \phi_k \d x ,  & \\ 
\label{eq:A matrix}
    \bfA|_{kj} &= \int_{\Ga_h}\!\!\! \nb_{\Ga_h} \phi_j \cdot \nb_{\Ga_h} \phi_k \d x
     &\ = &\ \sum_E \int_{E}\!\!\! \nb_{E} \phi_j \cdot \nb_{E} \phi_k \d x. &
\end{alignat}

\subsubsection{Computing local matrices}
\label{sec:local_matrix}
Using the integral formula \eqref{eq:basic integral formula} with \eqref{eq:basic function formula} we directly simplify the elementwise contributions of the local mass matrix as
\begin{align*}
    \bfM\loc|_{kj} :=&\ \int_E \phi_j \phi_k \d x =  \int_{\hat E} \vphi_{\hj} \vphi_{\hk} |\det DF_E| \d \xi  \qquad (\hj,\hk=1,2,3) .
\end{align*}

To tackle the gradient terms of the stiffness matrix we relate the gradients of the global and reference basis functions using \eqref{eq:basic function formula} and  \eqref{eq:discrete tangential gradient},
\begin{equation}
\label{eq:surf grad - pre}
	\nabla_E \phi_j (x) = \big(\id - \nu_E \nu_E^T\big) \nb \phi_j(x) = \big(\id - \nu_E \nu_E^T\big) DF_E^{-T} \nabla \vphi_{\hj}(\xi).
\end{equation}
The Jacobian $DF_E^{-T}$ is explicitly expressed utilising the element map \eqref{eq:F_map}, and reads:
\begin{align}
\label{eq:jacobian_inv}
	L_E^{-1} := DF_E^{-T} =\left[
	\begin{array}{c}
	 	\left(\sum_{\hi=1}^6 x_i\,\partial_x\phi_{\hat{i}}\right)^{T} \\[.3em]
	  	\left(\sum_{\hi=1}^6 x_i\,\partial_y\phi_{\hat{i}}\right)^{T} \\[.3em]
	 	\Bigl(\left(\sum_{\hi=1}^6 x_i\,\partial_x\phi_{\hat{i}}\right) \times \left(\sum_{\hi=1}^6 x_i\,\partial_y\phi_{\hat{i}}\right)\Bigr)^{T}
	\end{array}
	\right]^{-1} .
\end{align}
We further introduce the notation $C_E := L_E^{-1}$, and recall that $L_E \hat{\nu} = \nu_E$ (cf.~Figure~\ref{fig:notation}). 
Thus the projection for $\nabla_E \phi_j$ in \eqref{eq:surf grad - pre} is simplified as
\begin{align*}
	\nabla_E \phi_j (x) = C_E^T \nabla \vphi_{\hj}(\xi) - \nu_E \nu_E^T C_E^T \nabla \vphi_{\hj}(\xi) = C_E^T \nabla \vphi_{\hj}(\xi) - \nu_E \hat{\nu}^T \nabla \vphi_{\hj}(\xi) = C_E^T \nabla \vphi_{\hj}(\xi) , 
\end{align*}
where we used that $\hat{\nu}^T \nabla \vphi_{\hj} = 0$ for all reference basis functions. 

The formula for the local stiffness matrix then reads
\begin{equation}
\label{eq:Aloc - precise}
	\bfA\loc|_{kj} = \int_{E} \nb_{E} \phi_j \cdot \nb_{E} \phi_k \d x = \int_{\hat E} \nb \vphi_{\hj} \cdot C_E C_E^T \nb \vphi_{\hk} |\det DF_E| \d \xi .
\end{equation}

\subsubsection{Quadrature}
Since the polynomial degree of the integrands for the local mass matrix $\bfM\loc$ and for the stiffness matrix $\bfA\loc$ is $8$, we use high-order quadrature rules to compute these integrals (almost) exactly.  
Using a quadrature rule of order 8, with nodes and weights $(\xi_i,w_i)_{i=1}^Q$, we obtain:
\begin{align}
	\label{eq:reference_quadrature_mass}
	\bfM\loc|_{kj} = &\ \sum_{i=1}^Q w_i \vphi_{\hk}(\xi_i) \vphi_{\hj}(\xi_i) |\det L_E(\xi_i)| %+ \mathcal{O}(h^{q+1}) 
%	=: \sum_{i=1}^Q w_i  |\det L_E(\xi_i)| \bfM\rf|_{\hk \hj i}^q  %+ \mathcal{O}(h^{q+1})
	, \\ 
	\label{eq:reference_quadrature_stiffness}
	\bfA\loc|_{kj} = &\ \sum_{i=1}^Q w_i \nb \vphi_{\hj}(\xi_i) \cdot C_E C_E^T \nb \vphi_{\hk}(\xi_i) |\det L_E(\xi_i)| % + \mathcal{O}(h^{q+1}).
\end{align}

\section{Implementation}
\label{sec:implementation}
We first introduce some notation, $|\cdot|$ refers to the cardinality, i.e., the number of elements in a set. 
Let:
\begin{itemize}
\item $E$ be the set of all elements,
\item $N$ be the set of all nodes,
\item $d$ be the dimension of the element (and thus of the underlying manifold),
\item $p$ be the order of the local triangulation,
\item $\boldsymbol{N}$ be the coordinate array corresponding to the nodes $N$, 
	\par\hspace*{1em}$\Rightarrow$ $\dim (\mathbf{N}) = |N| \times \binom{p+d}{d}$,
\item and $\boldsymbol{E}$ be the connectivity array of the grid
	\par\hspace*{1em}$\Rightarrow$ $\dim (\mathbf{E}) = |E| \times \binom{p+d}{d}$. 
\end{itemize}

\subsection{Data structures}
\label{sec:data_structures}
The computational mesh is stored as the coordinate and connectivity arrays, listing each node from $N$ and respectively each element from $E$ (by listing the indices of nodes spanning the element). An example for a quadratic bulk mesh approximating a half disk is given in Figure~\ref{fig:2dp2_bulk} and Table~\ref{tab:elem_array}--\ref{tab:node_array}. 
Our implementation always assume that the nodes are listed anti-clockwise, first listing corners, then intermediate points, see Table~\ref{tab:elem_array}.
This notation is conform with usual finite element notation, see, e.g., \cite{Bartels_et_al_2006, funken2011, funken_beuter_2024}.
%For more insights we refer to \cite{Bartels_et_al_2006} who implemented curvilinear P2 elements for the flat case.

\begin{figure}[htbp]
  \centering
  % Left column (TikZ + Elements table stacked vertically)
  \begin{minipage}[t]{0.45\textwidth}
    \centering
    % TikZ picture
    \begin{tikzpicture}[scale=3,
	  every node/.style={
	    circle,
	    draw,
	    fill=gray!30,
	    minimum size=10pt,   % fixed node size
	    inner sep=0pt,       % prevent size variation by label length
	    font=\scriptsize
	  }]
	  
	  \pgfmathsetmacro{\a}{0.75}
	  
	  % Nodes
	  \coordinate (1) at (-1,0*\a);
	  \coordinate (2) at (0,0*\a);
	  \coordinate (3) at (0,1*\a);
	  \coordinate (4) at (-0.5,0*\a);
	  \coordinate (5) at (0,0.5*\a);
	  \coordinate (6) at (-0.707,0.707*\a);
	  \coordinate (7) at (1,0*\a);
	  \coordinate (8) at (0.5,0*\a);
	  \coordinate (9) at (0.707,0.707*\a);
%	  \coordinate (10) at (-0.5,0.25*\a);
%	  \coordinate (11) at (0,0.25*\a);
%	  \coordinate (12) at (0,0.75*\a);
%	  \coordinate (13) at (0.5,0.25*\a);
	
	  % Curves
	  \draw[thick, domain=-1:0, samples=50] plot (\x, {\a*(1 - \x - 2*\x*\x)});
	  \draw[thick, domain=0:1, samples=50] plot (\x, {\a*(1 + \x - 2*\x*\x)});
	  \draw[thick, gray, domain=-1:1, samples=100] plot (\x, {\a*sqrt(1 - \x*\x)});
	
	  % Lines
	  \draw[thick] (1) -- (2) -- (3);
	  \draw[thick] (2) -- (7);
%	  \draw[thick] (1) -- (5);
%	  \draw[thick] (5) -- (7);	
	
      \foreach \i in {1,...,9} \node at (\i) {\i}; 	
%      \foreach \i in {1,...,13} \node at (\i) {\i}; 			
	\end{tikzpicture}
    \caption{2D P2 bulk elements}
    \label{fig:2dp2_bulk}

    \vspace{1em}

    % Elements table below TikZ
    \begin{tabular}{cccccc}
      \toprule
      \multicolumn{6}{c}{$\mathbf{E}$ (Elements)} \\
      \midrule
      1 & 2 & 3 & 4 & 5 & 6 \\
      2 & 7 & 3 & 8 & 9 & 5 \\
%      1 & 2 & 5 & 4 & 11 & 10 \\
%      2 & 7 & 5 & 8 & 13 & 11 \\
%      5 & 7 & 3 & 13 & 9 & 12 \\
%      5 & 3 & 1 & 12 & 6 & 10 \\
      \bottomrule
    \end{tabular}
    \captionof{table}{Element array}
    \label{tab:elem_array}
  \end{minipage}
  \hfill
  % Right column (Nodes table taking whole height)
  \begin{minipage}[t]{0.45\textwidth}
    \centering
    \begin{tabular}{rr}
      \toprule
      \multicolumn{2}{c}{$\mathbf{N}$ (Nodes)} \\
      \midrule
      $-1$ & 0 \\
      0 & 0 \\
      0 & 1 \\
      $-0.5$ & 0 \\
      0 & 0.5 \\
      $-\sqrt{2}$ & $\sqrt{2}$ \\
      1 & 0 \\
      0.5 & 0 \\
      $-\sqrt{2}$ & $\sqrt{2}$ \\
%      $-0.5$ & 0.25 \\
%      0 & 0.25 \\
%      0 & 0.75 \\
%      0.5 & 0.25 \\
      \bottomrule
    \end{tabular}
    \captionof{table}{Coordinate array}
    \label{tab:node_array}
  \end{minipage}
\end{figure}

\subsection{Naive implementation, first optimizations}
%\label{sec:naiv_approach}
%The standard way to implement the assembly for finite elements, utilising the above setting is as follows: 
%\begin{enumerate}
%\item Preallocate the matrices and call precomputable quantities like the reference mass matrix, and evaluation of the product of reference basis functions on the quadrature nodes.
%\item Loop over all elements:
%\begin{enumerate}
%\item Determine map from local to global indices.
%\item Determine the Jacobian of $F$ given in \eqref{eq:F_map} and its determinant
%\item Compute local mass matrix by direct product of $|\det F|$, weight and precomputed reference mass
%\item Compute local stiffness matrix
%\begin{itemize}
%\item For each quadrature node, compute the inverse $C$ and determine the contribution of the quadrature node
%\item Sum and store elementwise result for local stiffness matrix
%\end{itemize}
%\end{enumerate}
%\item Allocate the mass and stiffness matrices based on the elementwise results.
%\end{enumerate} 

\label{sec:naive_approach}
For standard finite element implementations, the \emph{element-by-element} assembly of global mass and stiffness matrices is realised by looping over all elements, see, e.g., \cite{Braess_2007,Ern_Guermond_2004}. In each element, one obtains the local contributions as in Section~\ref{sec:local_matrix} and sums it into the correct place in the matrices. Typical computations are:
\begin{enumerate}
	\item Precompute reference quantities (e.g., basis function values and gradients at quadrature points).
	\item For each element, determine global indices, and the local element matrices.
	\item Accumulate the local contributions into the global sparse matrices.
\end{enumerate}

The structure of this setting is very clear and easy to implement, but the elementwise computation is inefficient in \matlab due to the use of loops (over elements and quadrature points), and more crucially repeated assignment to the global matrix is expensive. This prevents efficient compilation and usually leads to quadratic runtime $\mathcal{O}(|E|^2)$ in the number of elements $|E|$. 

There is a remedy to the elementwise approach, by storing the contributions of each element in cell arrays---storing the local mass and stiffness matrix, and index vectors (mapping local to global indices)---and assembling the global matrices only after computing all contributions, using the \matlab native \verb|sparse|. This allows to obtain almost linear runtime of $\mathcal{O}(|E|\log{|E|})$. See also the discussion in \cite[Section~3.3]{funken2011}. 
   
In general we can not expect to reduce the computation cost under $\mathcal{O}(|E| \log|E|)$ utilising \verb|sparse|. Main reason being that \matlab stores the sparse arrays in compressed sparse column format (for more details see, e.g., \cite{Cuvelier2013}) which requires duplicate free index arrays, which is not available with the above structure. The log-linear asymptotic arises solely from sorting and summing duplicate contributions.

High level FEM packages avoid this bottleneck by precomputing the sparsity structure and using better suited storage formats which allows for sort-free assignment of elementwise contributions, see, e.g.,~\cite{DUNE_article,ALBERTA_article,MFEM_article}.

\begin{remark}
Since for medium-scale applications ($|E| < 10^7$) the largest potential for optimization does not lie in dealing with other storage formats, cf.~Figure~\ref{fig:d2p2_surface}, we will continue to employ the coordinate list format and accept the small loss in time which is required for sorting into compressed sparse column format. 
\end{remark}

\subsection{Speed-up using \matlab's page-wise matrix functions}
In \matlab, the key to a faster runtime is vectorization via \matlab's \emph{page-wise matrix functions} \cite{croucher2024paged} first introduced in 2021. Allowing for parallel computation of all quantities which vastly decreases the computational cost. Through these functions the repetitive matrix operations of finite element assembly are implemented without any for-loops, instead we employ a highly efficient matrix--vector arithmetic: operations are automatically parallel, and allow for computations on a GPU by a simple command.

\begin{figure}
	\begin{center}
			\begin{tikzpicture}[every node/.style={anchor=north east,fill=white,minimum width=1.4cm,minimum height=7mm}]
				\matrix (mA) [draw,matrix of math nodes]
				{
				{|}&{|}&{|} \\
				x^{E_K}_1 & x^{E_K}_2 & x^{E_K}_3 \\
				{|}&{|}&{|} \\
				};
				\matrix (mB) [draw,matrix of math nodes] at ($(mA.south west)+(1.5,0.7)$)
				{
				{|}&{|}&{|} \\
				x^{E_1}_1 & x^{E_1}_2 & x^{E_1}_3\\
				{|}&{|}&{|} \\
				};
				\draw[dashed](mA.north east)--(mB.north east);
				\draw[dashed](mA.north west)--(mB.north west);
				\draw[dashed](mA.south east)--(mB.south east);
			\end{tikzpicture}
	\end{center}
	\caption{Coordinate array $\mathbf{N}$, where each page corresponds to some element $E_i$ for all elements ($i=1,\dots,|E|$). The vector $x_j^{E_i}$ stores all nodal coordinates of the $j$-th spatial dimension of the nodes of the given Element.}
	\label{fig:coordinate_array}
\end{figure}
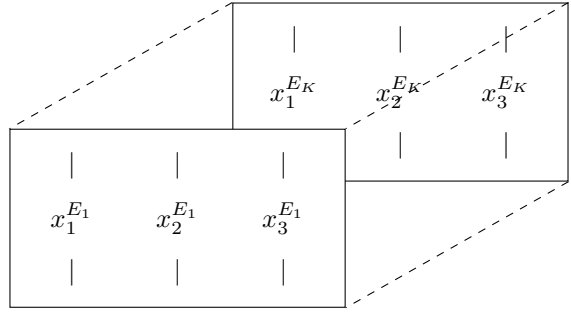

We will store all data in multidimensional arrays where each page corresponds to an element or a quadrature node, and perform all operations of the element-by-element matrix assembly via \emph{page-wise functions}.
As an example the coordinate array $\mathbf{N}$ which stores the coordinates for each node of each element for a linear 2D element would give the 3D array given in Figure~\ref{fig:coordinate_array}.

\subsubsection*{Main functions used for vectorization}
We have to carefully store and modify our data to realise the benefits of vectorization, in the following we briefly introduce the page-wise matrix functions used for data manipulations.

\begin{itemize}
	\item The $n$-dimensional array $A \in \mathbb{R}^{N_1 \times \dots \times N_n}$ is stored in \emph{column-major order}. As usual, the first two dimensions are labelled \emph{rows} and \emph{columns}, and the remaining dimensions are referred to as \emph{pages}. In the following, we typically deal with 3D arrays, and when quadrature is involved, we occasionally require 4D arrays.
	\item Multidimensional arrays can be directly assigned, for example:
	\[
	\mathtt{A(j,:,:)} = B, \quad B \in \mathbb{R}^{N_2 \times N_3},
	\]
	which fills the first row of each page of the array. Alternatively, the same array can be generated using
	\[
	\mathtt{C = repmat(B, r1, \dots, rn)},
	\]
	which constructs a new array $C$ consisting of $r_1, \dots, r_n$ copies of $B$ along the specified dimensions.
	
	\item We rearrange multidimensional arrays utilising
	\[
	\texttt{C = reshape(B,[r1,\dots,rn])} , \andquad  
	\texttt{C = permute(B,[n1, \dots, nm])}.
	\]
	The reshaping operator returns an array $C$ with the same values as in $B$ but the size of each dimension is given by $[r_1,\dots,r_n]$. The $\mathtt{permute}$ allows to reorder the dimensions by the provided permutations $[n_1, \dots, n_m]$.%, i.e.\ the prior dimension $n_1$ will correspond to the rows of the array $C$.
	\item Pagewise operations are the extension of standard matrix and vector operations, but work on each page separately, i.e.~acting on the first two dimensions. 
	
	Transposing, inverting and multiplying matrices is achieved, respectively, by
	\begin{alignat*}{3}
		& \texttt{C = pagetranspose(A)} & \qquad & C(:,:,i,j) = A(:,:,i,j)^T , \\
		& \texttt{C = pageinv(A)} & \qquad &  C(:,:,i,j) = A(:,:,i,j)^{-1} , \\
		& \texttt{C = pagemtimes(A,B)} & \qquad & C(:,:,i) = A(:,:,i,j) * B(:,:,i) .
	\end{alignat*}
%	For \texttt{pagemtimes} the dimensions of $A$ and $B$ must match from the third dimension onwards, or one of them has to be 1. 
%	%Clearly, if $A \in \mathbb{R}^{N_1 \times \dots \times N_n}$, then $B \in \mathbb{R}^{N_2 \times X \times \{1, N_3\} \times \dots}$.

%	First
%	\[ 
%	\texttt{C = pagetranspose(A)}
%	\]
%	transposes each page of $A$, i.e.\ for a 4D array we have $C(:,:,i,j) = A(:,:,i,j)^T$. Matrix and vector norms can be computed by
%	\[
%	\texttt{C = pagenorm(A,p)} \quad \text{and} \quad \texttt{C = vecnorm(A,p,d)}
%	\]
%	which returns the array $C$ of matrix or vector norms specified by the argument $p$ (note that for matrix norms this also includes the Frobenius norm, which is for example used for vectorizing the assembly for mean curvature flow in the numerical example of Section~\ref{sec:ex2_KLL}). For the vector norm we can specify along which dimension $d$ the norm is computed. 
%	To determine pagewise inverses we use
%	\[
%	\texttt{C = pageinv(A)},
%	\]
%	and finally
%	\[
%	\texttt{C = pagemtimes(A,B)},
%	\]
%	which computes the matrix product of $A(:,:,i,j) * B(:,:,i)$. From the third dimension onwards, the dimensions must match between $A$ and $B$ or one of them has to be 1. Clearly, if $A \in \mathbb{R}^{N_1 \times \dots \times N_n}$, then $B \in \mathbb{R}^{N_2 \times X \times \{1, N_3\} \times \dots}$.
	\item Matrix and vector norms can be computed by
	\[
	\texttt{C = pagenorm(A,p)} \quad \text{and} \quad \texttt{C = vecnorm(A,p,d)}
	\]
	which returns the array $C$ of matrix or vector norms (including the the Frobenius norm, 
%	which we use in vectorizing the assembly for mean curvature flow, 
	see Section~\ref{sec:ex2_KLL}). For the vector norm we can specify along which dimension $d$ the norm is computed. 

	\item Dimensions are collapsed (i.e.~dimension of length one are removed, for example, after summing all products of quadrature results and weights) by 
	\[
	\texttt{C = squeeze(A)}.
	\]
%	\item Dense vectors are created efficiently by
%	\[
%	\texttt{f = accumarray(v,w)},
%	\]
%	taking an index vector $v$ (with possible duplicates) and a value vector $w$, returning a vector $f$ which accumulates the possible duplicate indices and sum their values in $w$.
\end{itemize}

\subsubsection*{Reason and trade-off for speed-up}
The speeds-up of $\ell$FEM is based on the following ideas:
\begin{itemize}
	\item argument checking (for example when checking each page to apply the most efficient inverse builder, vs only checking it once); 
	\item better parallelization through page operations, than the automatic multi-thread routines for loops;
	\item optimized linear algebra for small matrices.
\end{itemize}

\begin{remark}[Batch computations]
\label{rem:trade_off}
Page functions are particularly efficient when working with a large number of elements at once. This, however, begs the question of available storage to work with large multi-dimensional dense arrays. In our experiments in Section~\ref{sec:num_ex} we observe that significant speed-up is possible on common hardware. 
Furthermore, large scale-problems can be tackled batchwise, see Section~\ref{sec:ex1_runtime}, achieving equivalent speed-ups. 
Due to the suboptimal runtime of $\ell$FEM, large-scale problems should be solved on dedicated FEM packages like DUNE, FeniCSx, MFEM, etc.
\end{remark}

\subsection{Loop-free implementation}
\label{section:detailed description}
We now give a detailed description of $\ell$FEM's loop-free fully vectorized isoparametric finite element assembly codes, using the two-dimensional second-order surface FEM as an example.
For readability some variables are abbreviated, e.g., $\verb|N| = \verb|Nodes|$ for the Nodes array, however, the provided codes usually use full name variables.

\begin{lstlisting}[caption={$\ell$FEM assembly code for quadratic isoparametric elements on a 2D surface},label={lst:code_ellFEM}]
	function [M,A] = assembly_surface_2D_P2(N,E)
	% (a) precomputations on reference element, and quadrature 
	p = 2; d = 2; ref = nchoosek(p+d,d); % dof on ref element
	[~,Mq, grad_Fq, W] = precompute_surface_2D_P2(); 
	Q = length(W); % number of quadrature nodes
	p1_Fq(1:ref,1:Q) = grad_Fq(1,:,:); %d_x \vphi(\xi_i)
	p2_Fq(1:ref,1:Q) = grad_Fq(2,:,:); %d_y \vphi(\xi_i)
	% (b) Index arrays and coordinate array
	nE = size(E,1); % number of elements
	I = repmat(E,1,ref)';
	J = repmat(permute(E,[3,2,1]),ref,1,1); 
	X = permute(reshape(N(E',:),[ref,nE,d+1]),[1,3,2]); 
	% (c) reference element transformations
	p1_Phi = pagemtimes(X,'transpose',reshape(p1_Fq,[ref,1,1,Q]),'none');
	p2_Phi = pagemtimes(X,'transpose',reshape(p2_Fq,[ref,1,1,Q]),'none'); 
	L3 = cross(p1_Phi,p2_Phi);
	det = squeeze(vecnorm(L3,2,1)); % surface elements
	L = [p1_Phi, p2_Phi, L3]; % Jacobians of element transformations
	% (d) mass matrix 
	Wdet = W .* det;
	MLoc = reshape(Mq,[ref^2 Q]) * Wdet'; % elementwise mass matrices
	M = sparse(I(:),J(:),MLoc(:)); % assigning mass matrix
	% (e) stiffness matrix
	C = pageinv(L); % inverting the Jacobians
	B = pagemtimes(C,'transpose', reshape(grad_Fq,[d+1,ref,1,Q]),'none'); 
	BtB = pagemtimes(B,'transpose',B,'none');
	ALoc = pagemtimes(permute(reshape(BtB,[ref^2 nE Q]),[1,3,2]), reshape(Wdet',[Q,1,nE])); % elementwise stifness matrices
	A = sparse(I(:),J(:),ALoc(:)); % assigning stifness matrix
\end{lstlisting}

%\begin{lstlisting}[caption={Assembly routine for a 2D surface with linear elements},label={code:surface_d2p1}]
%function [M,A,b] = surface_assembly_d2p1(N,E,f,deg_q)
%% Precomputable quantaties
%[M_loc, A_loc11, A_loc12, A_loc21, A_loc22, W, Fq] = precompute_d2p1(deg_quad); 
%nE = length(E);
%% Generate index vectors
%I = repmat(E,1,3)'; E2 = E';
%J = reshape(repmat(E2(:)',3,1),9,nE);
%X = permute(reshape(N(E2,:),[3,nE,3]),[1 3 2]);
%a = X(1,:,:); b = X(2,:,:); c = X(3,:,:);
%L = pagetranspose(reshape([b-a,c-a,cross(b-a,c-a)./vecnorm(cross(b-a,c-a))],[3,3,nE])); %Derivative of trafo
%detDF = L(1,1,:).*(L(2,2,:).*L(3,3,:) -L(2,3,:).*L(3,2,:))
%	   +L(1,2,:).*(L(2,3,:).*L(3,1,:)-L(2,1,:).*L(3,3,:))
%	   +L(1,3,:).*(L(2,1,:).*L(3,2,:)+L(2,2,:).*L(3,1,:)) 
%C = pageinv(L); %Inverse of L
%%% Mass matrix
%M_loc = detDF.*M_loc_0;
%M = sparse(I(:),J(:),M_loc(:));
%%% Stiffness matrix
%A_loc = detDF .* ((C(1,1,:).^2 + C(2,1,:).^2 + C(3,1,:).^2).*A_loc11 + (C(1,2,:).^2 + C(2,2,:).^2 + C(3,2,:).^2).*A_loc22 ...
%    + (C(1,1,:).*C(1,2,:) + C(2,1,:).*C(2,2,:) + C(3,1,:).*C(3,2,:)).*(A_loc12+A_loc21));
%A = sparse(I(:),J(:),A_loc(:));
%%% Load vector
%xq = pagemtimes(Fq,X); %coord of quadrature nodes for all E
%f_xq = reshape(f(xq(:,1,:),xq(:,2,:),xq(:,3,:)),[length(W),nE]);
%f_loc = pagemtimes(f_xq'.* (W.*squeeze(detDF)), F_eval); %quadrature
%b = accumarray(Elements(:),f_loc(:));
%\end{lstlisting}

\begin{itemize}
	\item[\textbf{(a)}] Precomputations of form functions and quadrature rule
	\begin{itemize}
	    \item[Line 2--5:] Precompute the product of all pairs of form functions evaluated at the quadrature nodes $(\xi_i)$ to obtain $\verb|Mq|$ for \eqref{eq:reference_quadrature_mass}, and a tensor containing evaluations of the gradients of the form functions $\verb|grad_Fq|=\nabla \vphi_{\hj} (\xi_i)$ required for $\bfA\loc$ \eqref{eq:reference_quadrature_stiffness}. Quadrature weights are stored in \verb|W|, and some constants are assigned: the order of the method $p$, the dimension of the elements $d$, and \verb|ref| degrees of freedom of the reference element.
		\item[Line 6--7:] The partial derivatives of the form functions at the quadrature nodes $\partial_x \vphi_{\hj}(\xi_i)$ are extracted from \verb|grad_Fq|, in preparation  for \eqref{eq:jacobian_inv}.
	\end{itemize}
	\item[\textbf{(b)}] Create index arrays $(I,J)$ and coordinate array $X = \mathbf{N}$.
	\begin{itemize}
	\item [Line 10--11:] The index arrays are directly constructed from the elements array $E$. 
%	This automatically constructs all index pairs of global coordinates for each element at once. Thus each consecutive set of 36 index pairs exactly matches all possible pairs of indices of one element. 
	These are used to map the reference indexing to the global indexing $(\hj, \hk) \leftrightarrow (j,k)$ for \eqref{eq:reference_quadrature_mass}--\eqref{eq:reference_quadrature_stiffness}.
	\item [Line 12:] Each page of $X$ corresponds to an element, cf.~Figure~\ref{fig:coordinate_array}.
	\end{itemize}
	\item[\textbf{(c)}] Reference element transformations and Jacobians
	\begin{itemize}
	\item [Line 14--18:] Compute each row and store the transformation matrices $L_E$ \eqref{eq:jacobian_inv}, calculate the surface elements via norm of the cross products.
	\end{itemize}
	\item[\textbf{(d)}] Mass matrix
	\begin{itemize}
	\item [Line 20--21:] Determine product of weights and surface elements, then a matrix product with \verb|Mq| sums each quadrature contribution to obtain $\bfM\loc$, cf.~\eqref{eq:reference_quadrature_mass}.
	\item [Line 22:] Assemble global mass matrix \eqref{eq:M matrix} utilising the index--value triple.
	\end{itemize}
	\item[\textbf{(e)}] Stiffness matrix
	\begin{itemize}
	\item [Line 24--25:] Determine the pagewise inverse and transpose of the Jacobians, cf.~\eqref{eq:jacobian_inv}. 
	\item [Line 26--27:] A paged products compute local stiffness matrix: sum each quadrature contribution and scale by weight to obtain the value of $\bfA\loc$, cf.~\eqref{eq:reference_quadrature_stiffness}.
	\item [Line 28:] Assemble global stiffness matrix \eqref{eq:A matrix} utilising the index--value triple.
	\end{itemize}
\end{itemize}

\section{Functions provided in the $\ell$FEM package}
\label{sec:overview_ellFEM}
The package $\ell$FEM provides codes for the assembly of the mass and stiffness matrices for the following bulk and surface domains, dimensions $d$, and elements of order $k = 1,2$ given in Table~\ref{tab:ellFEM_overview}.

\begin{table}[htbp]
	\centering
	\begin{tabular}{l c c c}
		\toprule
%		& \multicolumn{3}{c}{dimension} \\
%		\cmidrule(r){2-4}
		domain 				& 1D 		& 2D 		& 3D     \\ 
		\midrule
		surface 			& P1, P2 	& P1, P2 	&        \\
		bulk 	&  			& P1, P2 	& P1, P2 \\
		\bottomrule
	\end{tabular}
	\caption{The isoparametric finite elements supported by $\ell$FEM.}
	\label{tab:ellFEM_overview}
\end{table}

\subsection{Bulk and surface assembly codes}

On the git repository \href{http://go.upb.de/ellFEM}{\texttt{go.upb.de/ellFEM}} we provide the functions \verb|assembly_bulk_dD_Pk| and \verb|assembly_surface_dD_Pk| for bulk and surface domains, respectively (for various dimensions $d$ and orders $p$, see Table~\ref{tab:ellFEM_overview}). They take the two inputs \verb|Nodes| and \verb|Elements| as described in Section~\ref{sec:data_structures}. The computation follows the derivation of Section~\ref{sec:assembly}, and the description in Section~\ref{section:detailed description}.

%The 8 assembly functions are called by the ,
%\[
%\texttt{[M,A] = assembly\_type\_dD\_Pk(N,E)}.
%\]
%The data structures of both are described in Section~\ref{sec:data_structures}. 
%It returns the sparse global mass matrix $M$ based on \eqref{eq:M matrix}, and the global stiffness matrix $A$ based on \eqref{eq:A matrix}. The computation follows the derivation of Section~\ref{sec:assembly}.

\subsection{Precomputations}
Pre-computations are provided in each assembly code as subfunctions \verb|precompute_surface_dD_Pk|, \verb|precompute_bulk_dD_Pk|.

For P1 elements the precomputation codes return 
\[
\texttt{[M\_ref, A\_ref\_mn] = precompute\_type\_dD\_Pk()}.
\]
with the exact reference mass matrix and the building blocks of the stiffness matrices, i.e.~$\texttt{A\_ref\_mn(k,j)} = \int_{\hat{E}} \partial_m \vphi_{\hk} \partial_n \vphi_{\hj}$, with $m, n = 1, \dotsc, d$ and $\hk,\hj = 1,\dots, \binom{n+d}{d}$. %Thus there are $1+d^2$ outputs which are hardcoded and only called. 
Since for linear elements the element transformation is constant, using these matrices \eqref{eq:M matrix} and \eqref{eq:A matrix} are computed directly.

For P2 elements all integrals are computed using a sufficiently accurate quadrature rule: Gauss quadrature, Dunavant quadrature \cite{Dunavant1985}, and Jaskowiec and Sukumar quadrature \cite{JaskowiecSukumar}, of order 10, 8, and 7 in one, two, and three dimensions, respectively.

 To this end, the precomputation code
\[
\texttt{[F\_eval, M\_eval, grad\_F\_eval, W] = precompute\_type\_dD\_Pk()}
\] 
returns evaluations for form functions $\vphi_{\hk}(\xi_i)$ and their gradients, as well as the quadrature weights.
The quadrature rule is hardcoded, but it can easily be modified according to the users needs.
%Based on the quadrature and the respective form functions (see, e.g., \eqref{eq:form_functions_s2dp2}) and their gradients, it returns in order 
%\begin{itemize}
%\item an array storing the set of form functions evaluated at the quadrature nodes $\vphi_{\hk}(\xi_i)$,
%\item the three dimensional array $\bfM\rf^q|_{\hk \hj i}$,
%\item an array storing the set of gradients of form functions evaluated at the quadrature nodes $\partial_m \vphi_{\hk}(\xi_i)$, and
%\item the weight vector $W$.
%\end{itemize}
%The quadrature node indices range from $i = 1,\dots, Q$. 

\subsection{Tutorials and numerical experiments}

In addition to the main functions, the folder \verb|numerical experiments| provides all the numerical experiments from Section~\ref{sec:num_ex}.
Extensive examples which display the use-cases and complement this paper are given in \verb|tutorials|, to highlight a few: Dziuk's algorithm \cite{Dziuk90} and the KLL algorithm \cite{KLL2017} for mean curvature flow, a code for solving parabolic surface PDEs in a moving bulk domain, and codes for Poisson problems both for bulk and surface domains.

\subsection{Additional features and testing unit}
In the folder \verb|auxiliary| we provide some additional functions which are commonly used when working with surface and bulk finite elements. 

An approximation to the lift operator, the function
\[
\texttt{nodes = lift\_nD(x,d)},
\]
expects the coordinates of boundary nodes $x$ and a distance function $d$, and returns a coordinate array which corresponds to the lifted boundary nodes. Note that the lift is inexact but the error can be controlled by the hardcoded but modifiable tolerance $\texttt{TOL}$. This implementation is based on \cite{DistMesh}, see also \cite[Section~4.2]{Demlow_Dziuk_2007}, and can be used for surfaces of dimensions $n = 1,2,3$.

Another important tool is the generation of high-order polynomial meshes. For dimensions $d = 1,2,3$ and orders $p = 2,3,4$ we provide:
\[
\texttt{[N\_Pp, E\_Pp, E\_plot] = mesh\_preprocess\_dD(N, E, p)
},
\]
whose inputs are node and element arrays, \verb|N| and \verb|E|, respectively, of a simplicial mesh. Additionally it expects $\texttt{p}$, the target polynomial degree. The code returns in order: The node and element array corresponding to the $p$-th order mesh based on the input mesh. The output $\texttt{E\_plot}$ contains the elements of the corresponding P1 mesh, allowing for plotting via Matlab-native $\texttt{trisurf}$.

We provide a \emph{testing unit} in order to verify correctness of the code after installation or possible modifications for testing each case of Table~\ref{tab:ellFEM_overview}. This is achieved by computing errors with respect to a precomputed matrices using some meshes, which are also part of the testing unit.

\section{Numerical experiments}
\label{sec:num_ex}

\subsection{Runtime experiment for elliptic model problems}
\label{sec:ex1_runtime}
We perform runtime efficiency tests using elliptic model problems in bulk and surface domains. 
% Meshes were generated by consecutive red refinements combined with lifting. 
The following experiments report on assembly times against element count, and on efficiency in terms of errors against assembly times. 
All CPU tests are done on an Intel i7-1260P (32 GB RAM); GPU tests are done on a NVIDIA Titan RTX (24 GB RAM). Memory limits, which arises exactly due the trade-off discussed in Remark~\ref{rem:trade_off} are addressed via batch processing.

%for both the vectorized code and a suboptimal elementwise assembly (based on the direct approach given in ) and on the right we measure the efficiency of errors versus time spent of the vectorized implementation.

\subsubsection*{Elliptic surface PDEs in 2D with P2-elements}
\begin{figure}[htbp]
    \centering
    \includegraphics[width=1\linewidth]{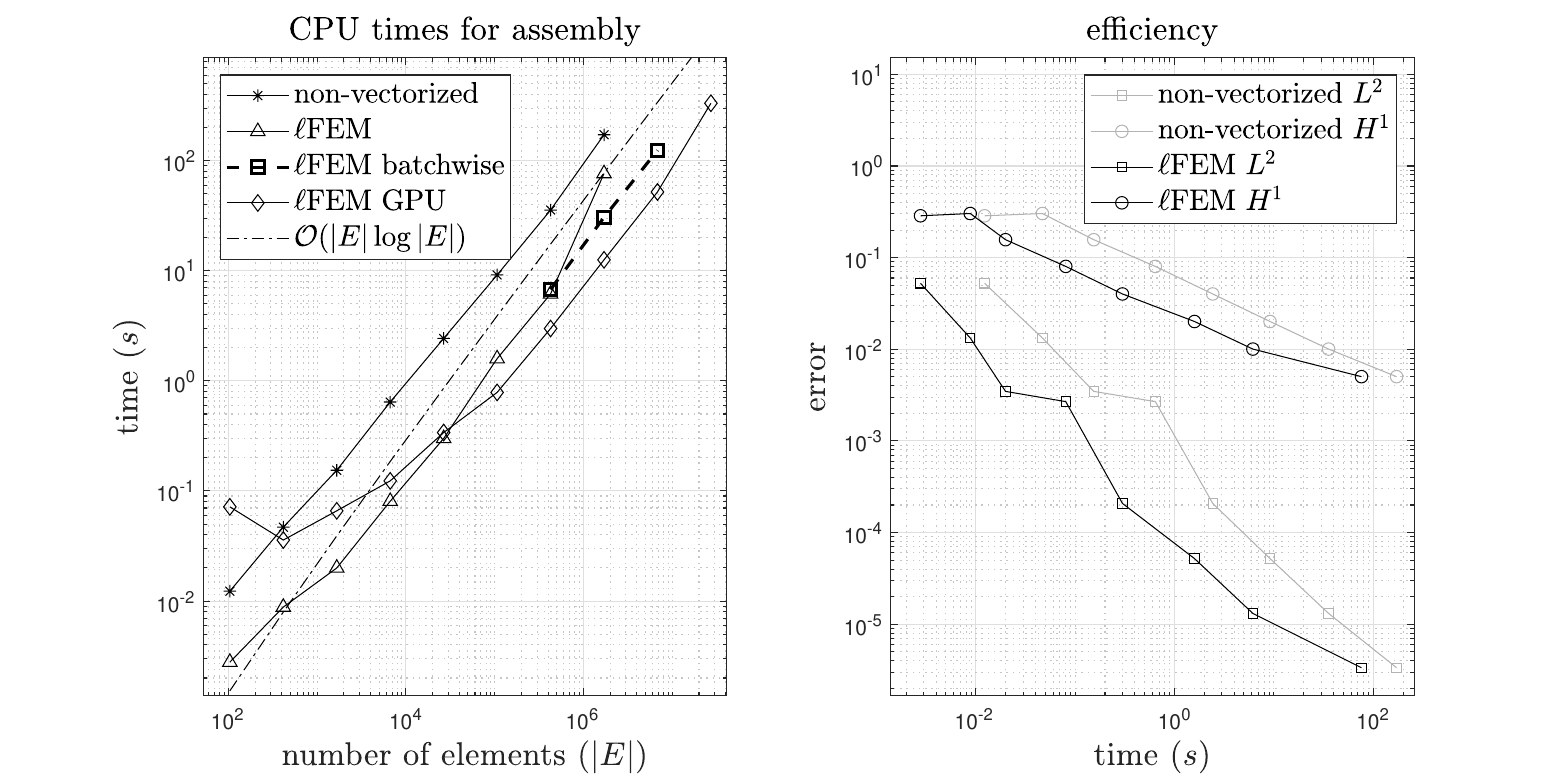}
    \caption{Assembly of mass, stiffness matrices and load vector for 2D surface FEM with P2 elements on (almost) uniform meshes using $\ell$FEM and non-vectorized codes. Computational cost (left), and time efficiency in $L^2$- and $H^1$-norm (right).}
	\label{fig:d2p2_surface}
\end{figure}

We solve the elliptic model problem \eqref{eq:model problem} on the unit sphere $\Ga = \mathbb{S}^2$ with the manufactured solution $u = x_1 x_2$ and inhomogeneity $f = 6 u$. 
The load vector is computed via Dunavant quadrature \cite{Dunavant1985} of order 8.
For uniqueness we seek the mean-value-free solution to the discrete formulation \eqref{eq:discretised model problem}. In Figure~\ref{fig:d2p2_surface} we compare $\ell$FEM ($\vartriangle$) and the non-vectorized code ($\ast$) (based on Section~\ref{sec:naive_approach}). The plot shows that the fully vectorized $\ell$FEM code is approximately one order of magnitude faster. 
RAM overflow beyond $\mathcal{O}(10^6)$ elements infers with the speed-up (see also Remark~\ref{rem:trade_off}), which is seen for the final data point of $\ell$FEM. Batchwise processing ($\square$), with batch size $M = 10^6$, resolves this issue.
We also tested GPU acceleration($\lozenge$) by providing the Nodes and Element arrays as $\texttt{gpuArray}$. Batchwise computation was required to avoid overflow on the GPU.
The GPU implementation outperforms all other methods only for high element counts, then a speed-up by a factor $2.5$ is observable. Since GPU computations are especially optimized for dense matrices, high-order methods and high-order quadratures are expected to benefit even more. 
We note that the batchwise computation is also limited by the overflow of storing the global matrices themselves.

\subsubsection*{Elliptic bulk PDEs in 2D with P1-elements}
\begin{figure}[htbp]
    \centering
    \includegraphics[width=1\linewidth]{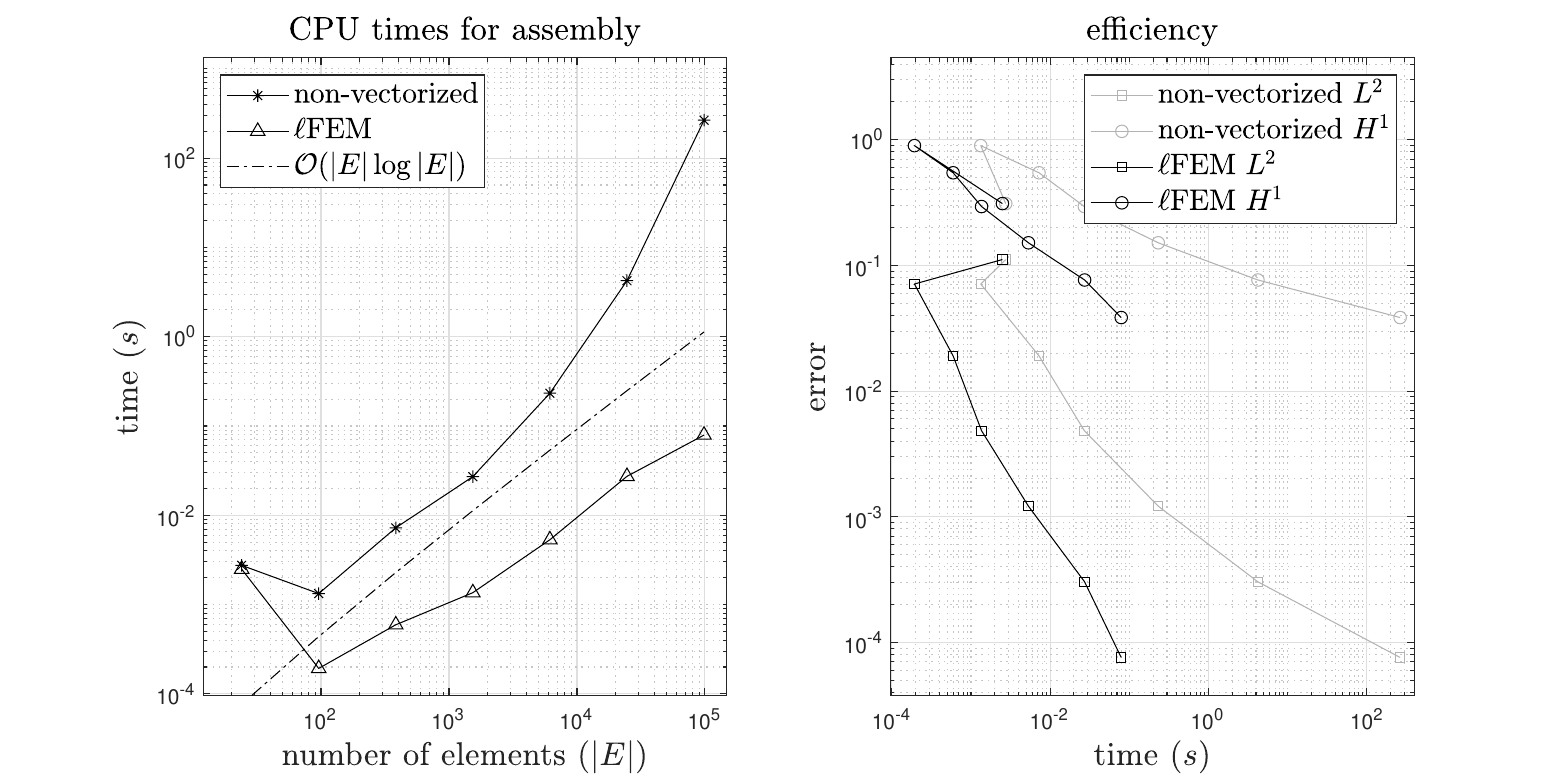}
    \caption{Assembly of mass, stiffness matrices and load vector for 2D bulk FEM with P1 elements on (almost) uniform meshes using $\ell$FEM and non-vectorized codes. Computational cost (left), and time efficiency in $L^2$- and $H^1$-norm (right).}
	\label{fig:d2p1_bulk}
\end{figure}

We solve a similar problem in the bulk $\Omega = \mathbb{B}^1$
\begin{align*}
-\Delta  u + \mu u = f \qquad \text{in} \ \Om , 
\end{align*}
with homogeneous Dirichlet boundary conditions. Assuming the radially symmetric solution $u = 1 - (x_1^2+x_2^2)^2$, we construct $f = 16 (x_1^2+x_2^2) + \mu u$. We choose $\mu = 10$ for the following numerical test.
%We only assembly the mass and stiffness matrix as a sufficient (first order) approximation of the load vector is given by $\bfM\bff$. 
Again we compare $\ell$FEM and the non-vectorized code. Figure~\ref{fig:d2p1_bulk} confirms that vectorizations yields one to two orders of magnitude speed-up. 

\subsection{Mean curvature flow: non-linear modifications and runtime comparison}
\label{sec:ex2_KLL}

We showcase the benefits and flexibility of $\ell$FEM via the \cite{KLL2017} algorithm for mean curvature flow. Note, that the finite element objects, including time-dependent mass and stiffness matrices, are required to be computed in each time step. Furthermore, we will demonstrate the implementation of non-linear expressions.

The algorithm for mean curvature flow developed in \cite{KLL2017} is based on the weak form of the coupled system:
\begin{equation}
\label{eq:mcf}
	\begin{alignedat}{3}
		v = &\ - H \nu, \qquad
		\partial^\bullet \nu &\ = &\ \Delta_{\Ga [X]} \nu + |\nabla_{\Ga [X]} \nu|^2 \nu,  \\
		\partial_t X  = &\ v \circ X ,\qquad
		\partial^\bullet H &\ = &\ \Delta_{\Ga [X]} H + |\nabla_{\Ga [X]} \nu|^2 H, 		
	\end{alignedat}
\end{equation}
where $X \colon \Ga^0 \times [0,T] \rightarrow \R^3$ denotes a parametrisation of $\Ga [X]$, $v$ is the normal velocity, and $H$ and $\nu$ are the mean curvature and outward unit normal, respectively.

Following \cite{KLL2017}, the evolving surface finite element / linearly implicit BDF discretisation leads to the following matrix--vector formulation, with $\bfx^n$ collecting the nodes of the discrete surface, and denoting $\bfu = (\bfn ; \bfH) \in \R^{4N}$:
\begin{equation*}
	\begin{aligned}
		\bfM(\widetilde{\bfx}^n)^{[4]} \dot{\bfu}^n + \bfA(\widetilde{\bfx}^n)^{[4]} \bfu^n = &\ \bff(\widetilde{\bfx}^n,\widetilde{\bfu}^n) \\
		\bfv^n = &\ - \bfH^n \bullet \bfn^n , \\
		\dot{\bfx}^n = &\ \bfv^n , 
	\end{aligned}
\end{equation*}
where $\dot{\bfx}^n$ denotes a BDF time derivative and $\widetilde{\bfx}^n$ the corresponding extrapolations. For complete details, including optimal-order fully discrete error estimates, we refer to \cite{KLL2017}. 

The non-linear term $\bff = (\bff_1 ; \bff_2)$ is determined by
\begin{align*}
	\bff_1(\bfx,\bfu)|_{j+\,( k -1)N} 
	= &\ \int_{\Ga_h[\bfx]} | \nabla_{\Ga_h[\bfx]} \nu_h |^2 \, (\nu_h)_k \, \phi_j[\bfx] , \\
	\bff_2(\bfx,\bfu)|_j 
	= &\ \int_{\Ga_h[\bfx]} | \nabla_{\Ga_h[\bfx]} \nu_h |^2 \, H_h \, \phi_j[\bfx],
\end{align*}
with $j = 1,\dotsc,N$ and $k = 1,2,3$.

We now demonstrate how the nonlinearities are implemented using $\ell$FEM. We focus on $\bff$:
\begin{lstlisting}
uLoc = permute(reshape(u(E,:).',[4 nE 6]),[3 1 2]);
alpha = squeeze(pagenorm(pagemtimes(C,'transpose', pagemtimes(reshape(grad_Fq,[3 6 1 16]), uLoc(:,1:3,:)),'none'),'fro'));
int = Wdet .* alpha.^2 .* permute( pagemtimes(Fq,uLoc),[3 1 2]);
fLoc = pagetranspose(pagemtimes(permute(int,[3 2 1]),Fq));
\end{lstlisting}
The implementation directly mirrors the theoretical formulations:
\begin{itemize}
	\item[Line 1:] Compute the pagewise variable \verb|uLoc| from the nodal values $\bfu = (\bfn,\bfH)$.
	\item[Line 2:] Compute $\alpha_h^2 := | \nabla_{\Ga_h[\bfx]} \nu_h |^2$, first evaluating the gradient $\nabla_{\Ga_h[\bfx]} \nu_h$ using precomputed quantities of Code~\ref{lst:code_ellFEM}, then computing its Frobenius norm using \verb|pagenorm|.
	%Integral transformation matrix $C$ multiplied with the gradient of the basis functions at the quadrature nodes and the nodal evaluations of the normal (as $\grad \nu = \bf\nu \grad \bf\phi$). 
	\item[Line 3:] Compute the pointwise product $\verb|int|_i = w_i \alpha_h^2(\xi_i) u_h(\xi_i)$.
	\item[Line 4:] Compute the pagewise matrix--vector product of $\verb|int|_i$ and $\phi_{\hk}(\xi_i)$ (which includes the summation of the quadrature rule).
\end{itemize}

We report on the MCF algorithm with P2 elements for a set of meshes approximating the initial surface $\mathbb{S}^2$, with final time $T=1$ and a fixed time step $\tau = 0.002$.
Table~\ref{fig:KLL_MCF_Assembly} reports on total times and on assembly contributions. We observe that the assembly times are mostly mesh-independent, for $\ell$FEM it takes slightly more than $60\%$ of the total computational cost, while for the non-vectorized implementation it is above $95\%$. The remaining computation time is mostly spent on solving the linear systems. In particular, for the finest mesh ($|E|=12008$) the total runtime decreased by a factor of more than $15$, and the assembly time by a factor of $25$.

\begin{table}[h]
\centering
\begin{tabular}{r S S S S}
\toprule
& \multicolumn{2}{c}{$\ell$FEM} & \multicolumn{2}{c}{non-vectorized} \\
\cmidrule(lr){2-3} \cmidrule(lr){4-5}
$|E|$ & {Total time [$s$]} & {Assembly [\%]} & {Total time [$s$]} & {Assembly [\%]} \\
\midrule
   212   & 0.398 & 72.6 & 2.975   & 96.3 \\
   440   & 0.746 & 68.0 & 5.875   & 95.9 \\
   956   & 2.001 & 71.5 & 14.678  & 96.1 \\
  1940   & 3.587 & 64.9 & 77.413  & 98.4 \\
  3728   & 7.431 & 67.3 & 122.287 & 98.0 \\
  6032   & 11.182 & 65.2 & 191.694 & 98.0 \\
  8192   & 21.214 & 62.6 & 306.851 & 97.4 \\
 10160   & 25.840 & 62.1 & 360.189 & 97.3 \\
 12008   & 26.738 & 61.2 & 421.007 & 97.5 \\
\bottomrule
\end{tabular}
\caption{Total and assembly time for the MCF algorithm \cite{KLL2017} using P2 elements.}
	\label{fig:KLL_MCF_Assembly}
\end{table}

\subsection{Dziuk's algorithm: MCF vs.~DUNE}
We briefly introduce Dziuk's algorithm \cite{Dziuk90} and compare our implementation, in particular the assembly runtime against the non-vectorized code, and an implementation in \href{https://www.dune-project.org}{DUNE}. Dziuk's algorithm \cite{Dziuk90} is the linearly implicit Euler / ESFEM discretisation of mean curvature flow, written as
\begin{equation*}
	\partial_t X(\cdot,t) = \laplace_{\Ga[X(\cdot,t)]} \textnormal{Id}_{\Ga[X(\cdot,t)]} ,
\end{equation*}
%\begin{align*}
%\int_{\Ga_h^n} \nabla_{\Ga_h^n} X_h^{n+1} \nabla_{\Ga_h^n} \phi_h + \frac{1}{\tau} X_h^{n+1} \phi_h = \int_{\Ga_h^n} \Big(\frac{1}{\tau} \textnormal{Id}_{\Ga_h^n} + 2 H \nu\Big) \phi_h.
%\end{align*}
%\todo{Notations?!}
%Both the test function $\phi_h$ and the discrete solutions $X_h^{n+1} \in S_h({\Ga_h^n})^3$.
%\bbk [kb: This is with spontaneous mean curvature $H$!! For MCF one has $H \equiv 0$.] \ebk 
which leads to the matrix--vector formulation (using the notation of Section~\ref{sec:ex2_KLL}):
\begin{equation*}
	\big(\bfM(\bfx^{n-1}) + \tau \bfA(\bfx^{n-1})\big) \bfx^n = \bfM(\bfx^{n-1}) \bfx^{n-1} .
\end{equation*}

The DUNE implementation is exactly the \href{https://www.dune-project.org/sphinx/content/sphinx/dune-fem/mcf_nb.html}{moving surface grid example in DUNE} (up to the linear solver, choosing $\theta = 1$ and enabling multi-threading). The initial mesh is given by collecting points on the sphere and mapping each point by $f(x,y,z) = (1 + 0.5\sin(2\pi(x+y))\cos(0.25z \pi))[x, y, z]^T$. We run the algorithms using P2 elements in the time interval $[0,0.3]$, with step size $\tau = 0.005$. We report on computation times for DUNE ($\square$), $\ell$FEM ($\vartriangle$), and the non-vectorized code ($\ast$).
We observe that DUNE is clearly faster than the non-vectorized code, but around $2$ to $4$ times slower than the $\ell$FEM implementation. We note that the improvements diminish for larger degrees of freedom due to the suboptimal asymptotic of $\ell$FEM.
\begin{figure}[htbp]
    \centering
    \includegraphics[width=1\linewidth]{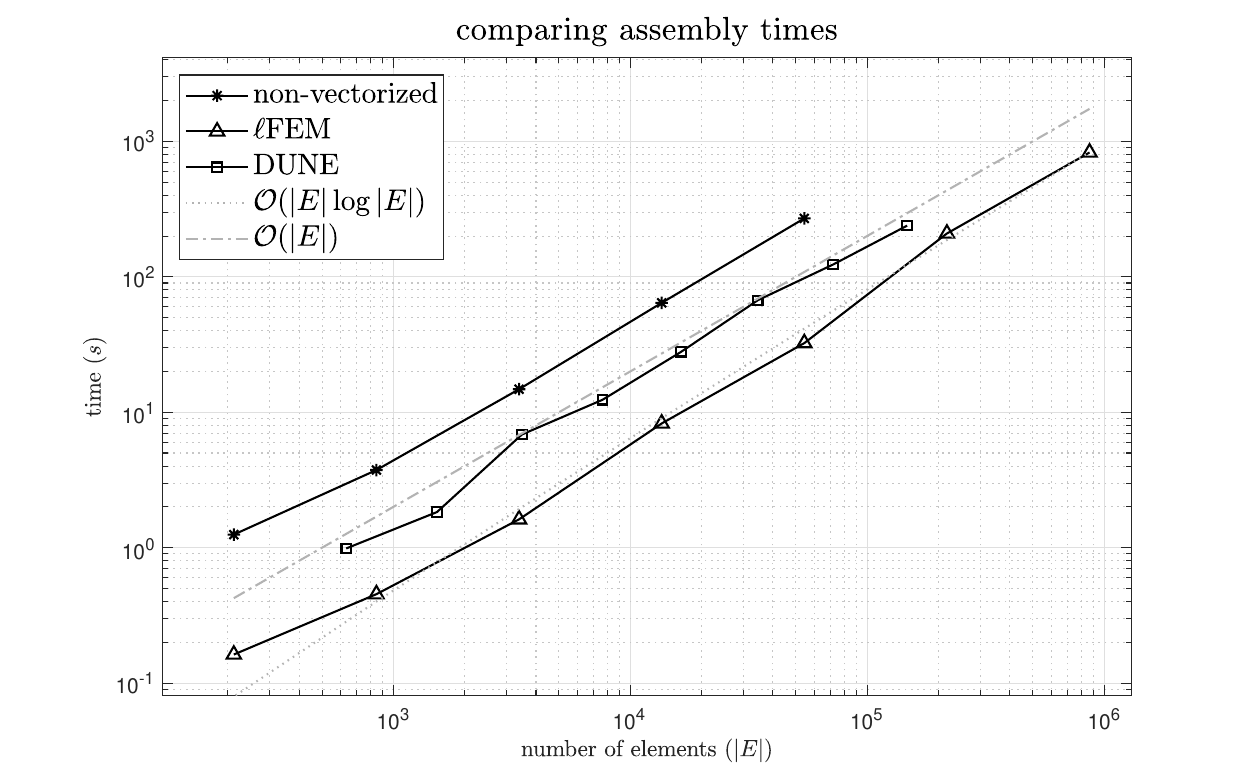}
    \caption{Comparing the assembly times for the non-vectorized, $\ell$FEM, and DUNE implementations of Dziuk's algorithm \cite{Dziuk90}. %The assembly times are cumulated over all time steps and plotted against the number of quadratic elements.
    }
	\label{fig:ex3}
\end{figure}

\subsection{Dimension analysis}
\begin{figure}
	\centering
	\begin{subfigure}[t]{0.48\textwidth}
		\centering
		\includegraphics[width=\linewidth]{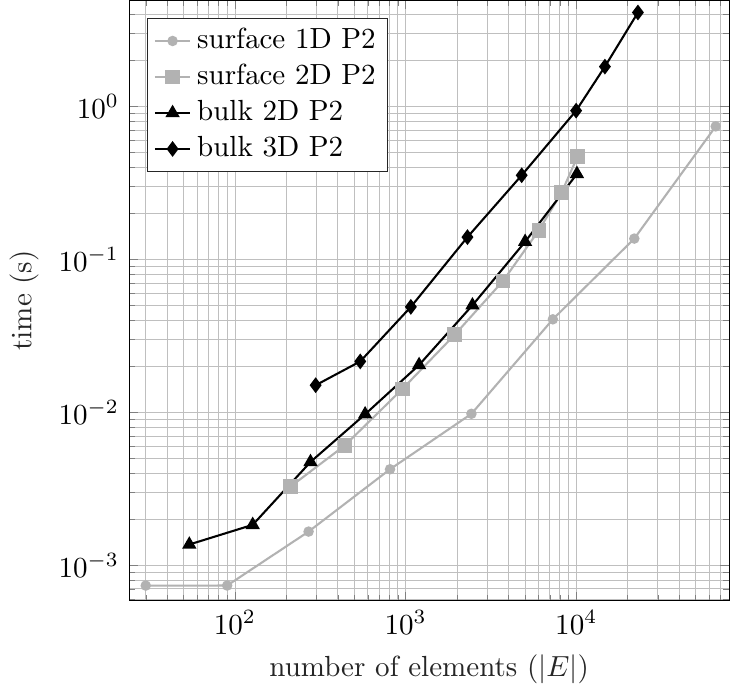}
		\caption{Runtime of assembly for meshes with varying number of elements $|E|$. Each line represents a quadratic finite element assembly of surface or bulk type, and the dimension given by the labels.}
		\label{fig:ex4_dim_order}
	\end{subfigure}
	\hfill
	\begin{subfigure}[t]{0.48\textwidth}
		\includegraphics[width=\linewidth]{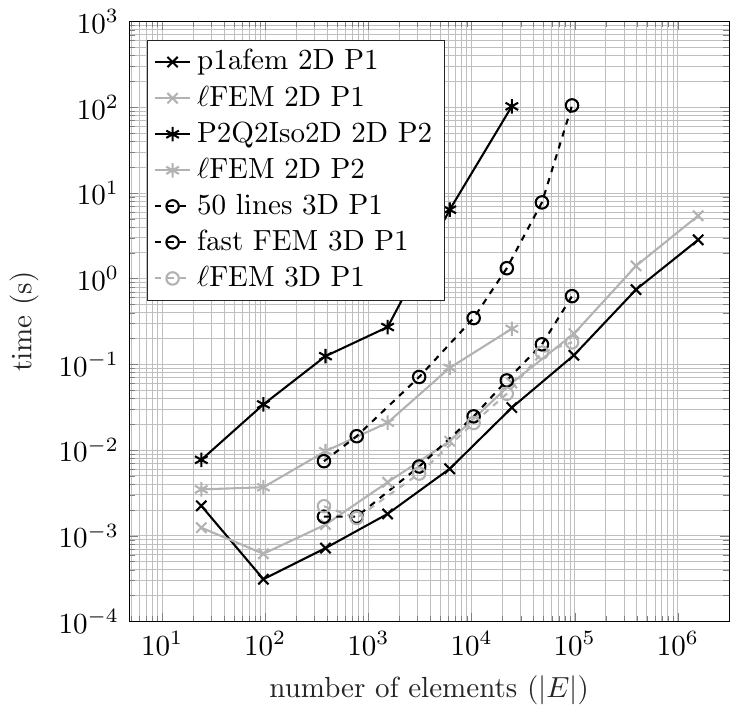} 
		\caption{Runtime of bulk assembly for various packages. The different order-- dimension pairs are marked by the same line marker. The $\ell$FEM runtimes are given in grey.}
		\label{fig:ex5_comparison}
	\end{subfigure}
\end{figure}

We now compare the P2 codes for dimensions $1$ to $3$, either for the unit ball (bulk), or for the unit sphere (surface) in the associated dimension. Figure~\ref{fig:ex4_dim_order} reports on the assembly times. The dimension-independence of the asymptotic runtime can be seen in the numerical results, and it also follows from analysing the provided algorithm. Note that the the runtimes for the two-dimensional bulk and surface match up, which demonstrates runtime independency from the underlying geometry.

%\begin{figure}
%	\centering
%	\begin{subfigure}[t]{0.48\textwidth}
%	    \centering
%	    \includegraphics[height=2.5in, trim={0 0 4.5in 0}, clip]{ex_4_5_test}
%	    \caption{Runtime of assembly for meshes with varying number of elements $\#E$. Each line represents a quadratic finite element assembly of surface or bulk type, and the dimension given by the labels.}
%		\label{fig:ex4_dim_order}
%	\end{subfigure}
%	\hfill
%	\begin{subfigure}[t]{0.48\textwidth}
%	    \centering
%	    \includegraphics[height=2.5in, trim={5in 0 0 0}, clip]{ex_4_5_test}   
%	    \caption{Runtime of bulk assembly for multiple packages. The different order-- dimension pairs are marked by the same line marker. The $\ell$FEM runtimes are given in grey.}
%		\label{fig:ex5_comparison}
%	\end{subfigure}
%\end{figure} 

\subsection{Comparison to other \matlab packages}
We perform a comparison experiment using established \matlab FEM packages. $\ell$FEM performs comparable -- often even better -- then these code, but retains the simple mathematical structure of the direct approach.

%We illustrate this by comparing the provided $\ell$FEM assembly codes against established finite element packages on a bulk domain. Note that we are able to compare to usual P1 finite element code as bulk P1 elements coincide with P1 elements as the discrete domain is identical). 
The following packages (giving their respective element order) are tested:
\begin{itemize}
	\item[-] \textbf{50 lines of MATLAB} \cite{Carstensen_1999}: P1 finite elements in 2D and 3D.
	\item[-] \textbf{P2Q2Iso2D} \cite{Bartels_et_al_2006}: Isoparametric P2 finite elements in 2D (allowing mixed triangle--quadrilateral meshes).
	\item[-] \textbf{p1afem} \cite{funken2011}: Efficient adaptive P1 finite element solver in 2D. Only the vectorized assembly is used here.
	\item[-] \textbf{fast FEM} \cite{Rahman_Valdman2013}: Vectorized assembly of isoparametric P1 elements in 2D and 3D.
\end{itemize}
Time comparison of these assembly process is in Figure~\ref{fig:ex5_comparison}.

Fast FEM \cite{Rahman_Valdman2013} and our code are comparably fast. Note that fast FEM implements vectorized functions, whereas $\ell$FEM utilises built-in \matlab functions (introduced after \cite{Rahman_Valdman2013} has appeared). The implementation p1afem \cite{funken2011} slightly outperforms $\ell$FEM for two reasons: It only computes the stiffness matrix. It uses a hardcoded optimization for constructing the entries of the stiffness matrix, which is slightly faster, but does not fully align with our development goals \textit{simplicity} and \textit{flexibility}, which are core goals for $\ell$FEM.

We are not aware of \matlab packages which efficiently implement surfaces finite elements in general, or efficient implementation of high-order bulk finite elements. As our codes are comparable in speed to other vectorized implementations of finite elements of low order, and the runtime order is independent of the order as seen in Figure~\ref{fig:ex4_dim_order}, we expect $\ell$FEM to perform comparatively in general.

\section*{Acknowledgement}

B.K.~would like to thank Dominik Edelmann for our long-standing and friendly competition on quick surface finite element implementation, and also for providing us a first version of our 3D mesh preprocessing tool. \redoff 

The authors would like to thank Andreas Dedner for his help regarding the fine-tuning of the DUNE implementation of Dziuk's algorithm.

The work of Bal\'azs Kov\'acs is funded by the Heisenberg Programme of the Deutsche Forschungsgemeinschaft (DFG, German Research Foundation) -- Project-ID 446431602,
and by the DFG Research Unit FOR 3013 \textit{Vector- and tensor-valued surface PDEs} (BA2268/6–1).

\bibliographystyle{siamplain}
\bibliography{assembly_literature}

@Article{ElliottRanner2021,
  author        = {C. M. Elliott and T. Ranner},
  journal       = {IMA Journal of Numerical Analysis},
  title         = {A unified theory for continuous-in-time evolving finite element space approximations to partial differential equations in evolving domains},
  year          = {2021},
  pages         = {1696--1845},
  volume        = {41},
}

@Article{EylesKingStyles2019,
  author     = {Eyles, J. and King, J. R. and Styles, V.},
  journal    = {Interfaces Free Bound.},
  title      = {A tractable mathematical model for tissue growth},
  year       = {2019},
  issn       = {1463-9963},
  number     = {4},
  pages      = {463--493},
  volume     = {21},
  doi        = {10.4171/ifb/428},
  fjournal   = {Interfaces and Free Boundaries. Mathematical Analysis, Computation and Applications},
  mrclass    = {35R35 (35Q92 65M60 74L15 92C30)},
  mrnumber   = {4046018},
}

@Article{StylesVanYperen2021,
  author        = {Styles, V. and Van Yperen, J.},
  journal       = {Num. Meth. PDE},
  title         = {Numerical analysis for a system coupling curve evolution attached orthogonally to a fixed boundary, to a reaction-diffusion equation on the curve},
  year          = {2021},
}

@Article{Bartels_et_al_2006,
  author  = {Bartels, S. and Carstensen, C. and Hecht, A.},
  journal = {J. Comput. Appl. Math.},
  title   = {{$\rm P2Q2Iso2D=2D$} isoparametric {FEM} in {M}atlab},
  year    = {2006},
  number  = {2},
  pages   = {219--250},
  volume  = {192},
}

@Article{DziukElliott_Acta_2013,
  author  = {Dziuk, G. and Elliott, C. M.},
  journal = {Acta Numerica},
  title   = {Finite element methods for surface {PDE}s},
  year    = {2013},
  pages   = {289--396.},
  volume  = {22},
  doi     = {10.1017/S0962492913000056},
}

@incollection {BGN_survey,
	AUTHOR = {Barrett, John W. and Garcke, Harald and N\"urnberg, Robert},
	TITLE = {Parametric finite element approximations of curvature-driven
	interface evolutions},
	BOOKTITLE = {Geometric partial differential equations. {P}art {I}},
	SERIES = {Handb. Numer. Anal.},
	VOLUME = {21},
	PAGES = {275--423},
	PUBLISHER = {Elsevier/North-Holland, Amsterdam},
	YEAR = {[2020] \copyright 2020},
	ISBN = {978-0-444-64003-1},
	MRCLASS = {65N30 (35K55 53Exx 74N05 76Txx 92C05)},
	MRNUMBER = {4378429},
}

@Article{DistMesh,
  author  = {Persson, P-O. and Strang, G.},
  journal = {SIAM Review},
  title   = {A Simple Mesh Generator in {MATLAB}},
  year    = {2004},
  number  = {2},
  pages   = {329--345.},
  volume  = {46},
}

@Article{DeckelnickDziukElliott_Acta_2005,
  author    = {Deckelnick, K. and Dziuk, G. and Elliott, C. M},
  journal   = {Acta Numerica},
  title     = {Computation of geometric partial differential equations and mean curvature flow},
  year      = {2005},
  pages     = {139--232},
  volume    = {14},
  publisher = {Cambridge University Press},
}

@Article{Demlow2009,
  author  = {A. Demlow},
  journal = {SIAM J. Numer. Anal.},
  title   = {Higher--order finite element methods and pointwise error estimates for elliptic problems on surfaces},
  year    = {2009},
  number  = {2},
  pages   = {805--807},
  volume  = {47},
  doi     = {10.1137/070708135},
}

@Article{Dunavant1985,
  author  = {Dunavant, D.A.},
  journal = {International journal for numerical methods in engineering},
  title   = {High degree efficient symmetrical {G}aussian quadrature rules for the triangle},
  year    = {1985},
  number  = {6},
  pages   = {1129--1148},
  volume  = {21},
}

@article {JaskowiecSukumar,
	AUTHOR = {Ja\'skowiec, Jan and Sukumar, N.},
	TITLE = {High-order symmetric cubature rules for tetrahedra and
	pyramids},
	JOURNAL = {Internat. J. Numer. Methods Engrg.},
	VOLUME = {122},
	YEAR = {2021},
	NUMBER = {1},
	PAGES = {148--171},
	ISSN = {0029-5981,1097-0207},
	MRCLASS = {65D30},
	MRNUMBER = {4196353},
	DOI = {10.1002/nme.6528},
}

@Article{Dziuk88,
  author  = {Dziuk, G.},
  journal = {Partial differential equations and calculus of variations, Lecture Notes in Math., 1357, Springer, Berlin},
  title   = {Finite elements for the {B}eltrami operator on arbitrary surfaces},
  year    = {1988},
  pages   = {142--155},
}

@Article{Dziuk90,
  author  = {Dziuk, G.},
  journal = {Numer. Math.},
  title   = {An algorithm for evolutionary surfaces},
  year    = {1990},
  number  = {1},
  pages   = {603--611},
  volume  = {58},
}

@Article{ElliottRanner2013,
  author  = {Elliott, C. M. and Ranner, T.},
  journal = {IMA J. Numer. Anal.},
  title   = {Finite element analysis for a coupled bulk-surface partial differential equation},
  year    = {2013},
  number  = {2},
  pages   = {377--402},
  volume  = {33},
}

@Article{MCF,
  author  = {Kov{\'a}cs, B. and Li, B. and Lubich, {Ch.}},
  journal = {Numer. Math.},
  title   = {A convergent evolving finite element algorithm for mean curvature flow of closed surfaces},
  year    = {2019},
  note    = {DOI:10.1007/s00211-019-01074-2},
  number  = {4},
  pages   = {797--853},
  volume  = {143},
}

@article{MCFdiff,
	title={Numerical analysis for the interaction of mean curvature flow and diffusion on closed surfaces},
	author={Elliott, C.~M. and Garcke, H. and Kov{\'a}cs, B.},
	journal={Numerische Mathematik}, 
	volume={151},
	number={4},
	pages={873--925},
	year={2022},
}

@article{CH_fulldiscrete,
	title={Error estimates for full discretization of {C}ahn--{H}illiard equation with dynamic boundary conditions},
	author={Bullerjahn, N. and Kov{\'a}cs, B.},
	journal={IMA J. Numer. Anal.},
	volume={46},
	number={2},
	pages={713--757},
	year={2026},
	doi={10.1093/imanum/draf009},
}

@article{adaptive_surface_FEM,
	title={A posteriori error estimates and space--time adaptivity for parabolic partial differential equations on stationary surfaces},
	author={Kov{\'a}cs, B. and Lantelme, M.},
	journal={to appear in SIAM J. Numer. Anal.},
	volume={},
	number={},
	pages={},
	year={2026},
	note={arXiv:2407.02101},
}

@article{bulksurface_coupling,
	title={Numerical analysis of an evolving bulk--surface model of tumour growth},
	author={Edelmann, D. and Kov{\'a}cs, B. and Lubich, {Ch.}},
	journal={IMA Journal of Numerical Analysis},
	volume={},
	number={},
	pages={},
	year={2024},
	doi={10.1093/imanum/drae077},
}

@article {funken2011,
    AUTHOR = {Funken, Stefan and Praetorius, Dirk and Wissgott, Philipp},
     TITLE = {Efficient implementation of adaptive {P}1-{FEM} in {M}atlab},
   JOURNAL = {Comput. Methods Appl. Math.},
  FJOURNAL = {Computational Methods in Applied Mathematics},
    VOLUME = {11},
      YEAR = {2011},
    NUMBER = {4},
     PAGES = {460--490},
      ISSN = {1609-4840,1609-9389},
   MRCLASS = {65N30 (65N50 65Y15)},
  MRNUMBER = {2875100},
       DOI = {10.2478/cmam-2011-0026},
}

@online{croucher2024paged,
  author       = {Mike Croucher},
  title        = {Paged Matrix Functions in MATLAB (2024 Edition)},
  year         = {2024},
  month        = {May 29},
  url          = {https://blogs.mathworks.com/matlab/2024/05/29/paged-matrix-functions-in-matlab-2024-edition/},
  note         = {MATLAB \& Simulink - MathWorks Blogs}
}

@article{Cuvelier2013,
author = {Cuvelier, François and Japhet, Caroline and Scarella, Gilles},
year = {2013},
month = {05},
pages = {},
title = {An efficient way to perform the assembly of finite element matrices in Matlab and Octave}
}

@article {KLL2017,
    AUTHOR = {Kov\'acs, Bal\'azs and Li, Buyang and Lubich, Christian},
     TITLE = {A convergent evolving finite element algorithm for mean
              curvature flow of closed surfaces},
   JOURNAL = {Numer. Math.},
  FJOURNAL = {Numerische Mathematik},
    VOLUME = {143},
      YEAR = {2019},
    NUMBER = {4},
     PAGES = {797--853},
      ISSN = {0029-599X,0945-3245},
   MRCLASS = {65M60 (35R01 65M12 65M15)},
  MRNUMBER = {4026373},
       DOI = {10.1007/s00211-019-01074-2},
}

@article {Rahman_Valdman2013,
    AUTHOR = {Rahman, Talal and Valdman, Jan},
     TITLE = {Fast {MATLAB} assembly of {FEM} matrices in 2{D} and 3{D}:
              nodal elements},
   JOURNAL = {Appl. Math. Comput.},
  FJOURNAL = {Applied Mathematics and Computation},
    VOLUME = {219},
      YEAR = {2013},
    NUMBER = {13},
     PAGES = {7151--7158},
      ISSN = {0096-3003,1873-5649},
   MRCLASS = {65N30},
  MRNUMBER = {3030557},
       DOI = {10.1016/j.amc.2011.08.043},
}

@article {funken_beuter_2024,
    AUTHOR = {Beuter, Stefanie and Funken, Stefan A.},
     TITLE = {Efficient {P}1-{FEM} for any space dimension in {M}atlab},
   JOURNAL = {Comput. Methods Appl. Math.},
  FJOURNAL = {Computational Methods in Applied Mathematics},
    VOLUME = {24},
      YEAR = {2024},
    NUMBER = {2},
     PAGES = {283--315},
      ISSN = {1609-4840,1609-9389},
   MRCLASS = {65M60},
  MRNUMBER = {4726526},
       DOI = {10.1515/cmam-2022-0239},
}

@misc{mallesham2025,
      title={Vectorized 3D mesh refinement and implementation of primal hybrid FEM in MATLAB}, 
      author={Harish Nagula Mallesham and Sharat Gaddam and Jan Valdman and Sanjib Kumar Acharya},
      year={2025},
      eprint={2509.11133},
      archivePrefix={arXiv},
      primaryClass={math.NA},
}

@article{anisotropicMCF,
	title={Error estimates for a surface finite element method for anisotropic mean	curvature flow},
	author={Deckelnick, K. and Garcke, H. and Kov{\'a}cs, B.},
	journal={}, 
	volume={},
	number={},
	pages={},
	year={2025},
	note={arxiv:2508.01078},
}

@article{numerical_surgery,
	title={Numerical surgery for mean curvature flow of surfaces},
	author={Kov{\'a}cs, B.},
	journal={SIAM Journal on Scientific Computing}, 
	volume = {46},
	number = {2},
	pages = {A645--A669},
	year = {2024},
	doi={10.1137/22M1531919},
}

@article{Demlow_Dziuk_2007,
	author = {Demlow, Alan and  Dziuk, Gerhard},
	title = {An Adaptive Finite Element Method for the Laplace–Beltrami Operator on Implicitly Defined Surfaces},
	journal = {SIAM Journal on Numerical Analysis},
	volume = {45},
	number = {1},
	pages = {421-442},
	year = {2007},
	doi={10.1137/050642873},
}

@article {Carstensen_1999,
    AUTHOR = {Alberty, Jochen and Carstensen, Carsten and Funken, Stefan A.},
     TITLE = {Remarks around 50 lines of {M}atlab: short finite element implementation},
   JOURNAL = {Numer. Algorithms},
  FJOURNAL = {Numerical Algorithms},
    VOLUME = {20},
      YEAR = {1999},
    NUMBER = {2-3},
     PAGES = {117--137},
      ISSN = {1017-1398,1572-9265},
   MRCLASS = {65N30 (65M60)},
  MRNUMBER = {1709562},
       DOI = {10.1023/A:1019155918070},
}

@techreport{Chen_2009_ifem,
  author       = {Long Chen},
  title        = {iFEM: an Integrated Finite Element Methods Package in {MATLAB}},
  institution  = {University of California at Irvine},
  year         = {2009},
  url          = {https://github.com/lyc102/ifem},
  note         = {Technical Report}
}

@article{Feifel_Funken_2024,
title = {Efficient P1-FEM for Any Space Dimension in Matlab},
author = {Stefanie Feifel and Stefan A. Funken},
pages = {283--324},
volume = {24},
number = {2},
journal = {Computational Methods in Applied Mathematics},
doi = {10.1515/cmam-2022-0239},
year = {2024},
lastchecked = {2025-10-22}
}

@book{Braess_2007, 
place={Cambridge},
edition={3},
title={Finite Elements: Theory, Fast Solvers, and Applications in Solid Mechanics},
publisher={Cambridge University Press},
author={Braess, Dietrich},
year={2007}
}

@book{Ern_Guermond_2004,
  address = {New York},
  author = {Ern, Alexandre and Guermond, Jean-Luc},
  doi = {10.1007/978-1-4757-4355-5},
  isbn = {978-0-387-20574-8},
  publisher = {Springer},
  series = {Applied Mathematical Sciences},
  title = {Theory and Practice of Finite Elements},
  volume = 159,
  year = 2004
}

@article {Bonito_Devore_Nochetto_2013,
    AUTHOR = {Bonito, Andrea and DeVore, Ronald A. and Nochetto, Ricardo H.},
     TITLE = {Adaptive finite element methods for elliptic problems with
              discontinuous coefficients},
   JOURNAL = {SIAM J. Numer. Anal.},
  FJOURNAL = {SIAM Journal on Numerical Analysis},
    VOLUME = {51},
      YEAR = {2013},
    NUMBER = {6},
     PAGES = {3106--3134},
      ISSN = {0036-1429,1095-7170},
   MRCLASS = {65N30 (65N50)},
  MRNUMBER = {3129757},
MRREVIEWER = {V\'it\ Dolej\v s\'i},
       DOI = {10.1137/130905757},
}

@article{DUNE_article,
author = {Bastian, Peter and Blatt, Markus and Dedner, Andreas and Engwer, Christian and Klöfkorn, R. and Kornhuber, Ralf and Ohlberger, Mario and Sander, Oliver},
year = {2008},
month = {01},
pages = {121-138},
title = {A generic grid interface for parallel and adaptive scientific computing. Part II: Implementation and tests in DUNE},
volume = {82},
journal = {Computing},
doi = {10.1007/s00607-008-0004-9}
}

@article{ALBERTA_article,
author = {Schmidt, Alfred and Siebert, Kunibert},
year = {2005},
month = {01},
pages = {},
title = {Design of adaptive finite element software: The finite element toolbox ALBERTA}
}

@article{MFEM_article,
title = {MFEM: A modular finite element methods library},
journal = {Computers \& Mathematics with Applications},
volume = {81},
pages = {42-74},
year = {2021},
note = {Development and Application of Open-source Software for Problems with Numerical PDEs},
issn = {0898-1221},
doi = {https://doi.org/10.1016/j.camwa.2020.06.009},
author = {Robert Anderson and Julian Andrej and Andrew Barker and Jamie Bramwell and Jean-Sylvain Camier and Jakub Cerveny and Veselin Dobrev and Yohann Dudouit and Aaron Fisher and Tzanio Kolev and Will Pazner and Mark Stowell and Vladimir Tomov and Ido Akkerman and Johann Dahm and David Medina and Stefano Zampini},
}

@article {dynbc,
	AUTHOR = {Kov\'acs, Bal\'azs and Lubich, Christian},
	TITLE = {Numerical analysis of parabolic problems with dynamic boundary
	conditions},
	JOURNAL = {IMA J. Numer. Anal.},
	FJOURNAL = {IMA Journal of Numerical Analysis},
	VOLUME = {37},
	YEAR = {2017},
	NUMBER = {1},
	PAGES = {1--39},
	ISSN = {0272-4979,1464-3642},
	MRCLASS = {65M60 (65M15)},
	MRNUMBER = {3614879},
	MRREVIEWER = {Jan\ Giesselmann},
	DOI = {10.1093/imanum/drw015},
}

@article {HochbruckHippStoher,
	AUTHOR = {Hipp, David and Hochbruck, Marlis and Stohrer, Christian},
	TITLE = {Unified error analysis for nonconforming space discretizations
	of wave-type equations},
	JOURNAL = {IMA J. Numer. Anal.},
	FJOURNAL = {IMA Journal of Numerical Analysis},
	VOLUME = {39},
	YEAR = {2019},
	NUMBER = {3},
	PAGES = {1206--1245},
	ISSN = {0272-4979,1464-3642},
	MRCLASS = {65M60 (65M12 65M15)},
	MRNUMBER = {3984056},
	MRREVIEWER = {Dami\'an\ P.\ Ginestar},
	DOI = {10.1093/imanum/dry036},
}
\end{document}